%% file: badia_li_martin_2025.tex
\documentclass[oneside,reqno,11pt,a4paper]{amsart}
\usepackage[T1]{fontenc}
\usepackage[margin=2.4cm]{geometry}
\usepackage{amsmath,amsfonts,amsbsy,amssymb,amscd,amsthm}
\usepackage[usenames,dvipsnames,svgnames]{xcolor}
\usepackage[utf8]{inputenc}
\usepackage[ruled]{algorithm2e}
\usepackage{graphicx}
\usepackage{svg}
\usepackage{caption}
\usepackage{subcaption}
\usepackage{doi}
\usepackage[foot]{amsaddr}
\usepackage{makecell}
\usepackage{hyperref}
\usepackage{acronym}
\usepackage{listings}
\usepackage{textgreek}
\usepackage{url}
\usepackage{color}
\usepackage{framed}
\usepackage{verbatim}
\usepackage{fancyvrb}
\usepackage{stmaryrd}
\usepackage{newtxtext}
\usepackage{newtxmath}
\usepackage{microtype} 
\usepackage[final]{pdfpages}
\usepackage{tikz}
\usepackage{booktabs}
\usepackage{multirow}

\usepackage[backend=biber,
sorting=none,
defernumbers=false,
maxbibnames=5,
style=numeric-comp,
isbn=false,
bibencoding=utf8,
safeinputenc, 
url=false,
doi=true,
eprint=false,
giveninits=true]{biblatex}
\hypersetup{
breaklinks=true,
bookmarksopen=true,
pdftitle={Unfitted Finite Element Interpolated Neural Networks},    
pdfauthor={Wei Li, Alberto F. Mart{\'{\i}}n, and Santiago Badia},     
colorlinks=true,       
linkcolor=blue,          
citecolor=blue,        
filecolor=black,      
urlcolor=blue           
}  
\definecolor{plum}{RGB}{160,43,147}
\definecolor{lightgreen}{RGB}{140,181,145}
\definecolor{lightblue}{RGB}{220,234,247}
\definecolor{lightred}{RGB}{251,227,214}
  
\tikzset{
  cross/.pic = {
  \draw[rotate = 45,thick] (-#1,0) -- (#1,0);
  \draw[rotate = 45,thick] (0,-#1) -- (0, #1);
  }
}

\newtheorem{theorem}{Theorem}[section]

\newtheorem{remark}[theorem]{Remark}

\acrodef{pde}[PDE]{partial differential equation}
\acrodef{fe}[FE]{finite element}
\acrodef{fem}[FEM]{finite element method}
\acrodef{agfem}[AgFEM]{aggregated \ac{fem}}
\acrodef{fdm}[FDM]{finite difference method}
\acrodef{fvm}[FVM]{finite volume method}
\acrodef{dof}[DoF]{degree of freedom}
\acrodefplural{dof}[DoFs]{degrees of freedom}
\acrodef{nn}[NN]{neural network}
\acrodef{cnn}[CNN]{convolutional \ac{nn}}
\acrodef{pinn}[PINN]{physics-informed \ac{nn}}
\acrodef{vpinn}[VPINN]{variational \ac{pinn}}
\acrodef{ivpinn}[IVPINN]{interpolated \ac{vpinn}}
\acrodef{feinn}[FEINN]{\ac{fe} interpolated \ac{nn}}
\acrodef{pfnn}[PFNN]{penalty-free \ac{nn}}
\acrodef{ibn}[IBN]{irregular boundary network}
\acrodef{gmg}[GMG]{geometric multigrid}
\acrodef{spd}[SPD]{symmetric positive definite}
\acrodef{stl}[STL]{stereolithography}
\acrodef{cae}[CAE]{computer-aided engineering}
\acrodef{hidenn}[HiDeNN]{hierarchical deep-learning \ac{nn}}
\acrodef{fenni}[FENNI]{\ac{fe} \ac{nn} interpolation}
\acrodef{pgd}[PGD]{proper generalised decomposition}

\graphicspath{{../plots}{./}}
\newcommand{\fig}[1]{Fig.~\ref{#1}}
\newcommand{\tab}[1]{Tab.~\ref{#1}}
\newcommand{\sect}[1]{Sect.~\ref{#1}}

\newcommand{\norm}[1]{\left\lVert #1 \right\rVert}
\newcommand{\ltwonorm}[1]{\left\lVert #1 \right\rVert _{L^2(\Omega)}}
\newcommand{\honenorm}[1]{\left\lVert #1 \right\rVert _{H^1(\Omega)}}
\newcommand{\jump}[1]{\lBrack #1 \rBrack}
\newcommand{\normal}{\pmb{n}}

\newcommand{\argmin}[1]{\underset{#1}{\mathrm{arg\,min}}\,}

\begin{document}

\title[Unfitted finite element interpolated neural networks]{Unfitted finite element interpolated neural networks}
\author{Wei Li$^1$}
\email{wei.li@monash.edu}
\author{Alberto F. Mart\'{\i}n$^2$}
\email{alberto.f.martin@anu.edu.au}
\author{Santiago Badia$^{1,*}$}
\email{santiago.badia@monash.edu}
\address{$^1$ School of Mathematics, Monash University, Clayton, Victoria 3800, Australia.}
\address{$^2$ School of Computing, The Australian National University, Canberra ACT 2600, Australia.}
\address{$^*$ Corresponding author.}

\date{\today}

\begin{abstract}
  We present a novel approach that integrates unfitted finite element methods and neural networks to approximate partial differential equations on complex geometries.
  Easy-to-generate background meshes (e.g., a simple Cartesian mesh) that cut the domain boundary (i.e., they do not conform to it) are used to build suitable trial and test finite element spaces.
  The method seeks a neural network that, when interpolated onto the trial space, minimises a discrete norm of the weak residual functional on the test space associated to the equation.
  As with unfitted finite elements, essential boundary conditions are weakly imposed by Nitsche's method. The method is robust to variations in Nitsche coefficient values, and to small cut cells.
  We experimentally demonstrate the method's effectiveness in solving both forward and inverse problems across various 2D and 3D complex geometries, including those defined by implicit level-set functions and explicit stereolithography meshes.
  For forward problems with smooth analytical solutions, the trained neural networks achieve several orders of magnitude smaller $H^1$ errors compared to their interpolation counterparts. These interpolations also maintain expected $h$- and $p$-convergence rates.
  Using the same amount of training points, the method is faster than standard PINNs (on both GPU and CPU architectures)  while achieving similar or superior accuracy.
  Moreover, using a discrete dual norm of the residual (achieved by cut cell stabilisation) remarkably accelerates neural network training and further enhances robustness to the choice of Nitsche coefficient values.
  The experiments also show the method's high accuracy and reliability in solving inverse problems, even with incomplete observations.
\end{abstract}


\maketitle

\input{intro.tex}
\input{method.tex}
\input{implmt.tex}
\input{numexp.tex}

\section*{Acknowledgments}
This research was partially funded by the Australian Government through the Australian Research Council (project numbers DP210103092 and DP220103160). This work was also supported by computational resources provided by the Australian Government through NCI under the NCMAS and ANU Merit Allocation Schemes.
W. Li gratefully acknowledges the Monash Graduate Scholarship from Monash University, the NCI computing resources provided by Monash eResearch through Monash NCI scheme for HPC services, and the support from the Laboratory for Turbulence Research in Aerospace and Combustion (LTRAC) at Monash University through the use of their HPC Clusters.

\section*{Declaration of generative AI and AI-assisted technologies in the writing process}
During the preparation of this work the authors used ChatGPT in order to improve language and readability. After using this tool/service, the authors reviewed and edited the content as needed and take full responsibility for the content of the publication.

\printbibliography

\end{document}

%% file: intro.tex
\section{Introduction} \label{sec:intro}
\Acp{pde} are powerful mathematical tools for modelling a wide range of physical phenomena, including heat transfer, fluid dynamics, and electromagnetics. Since obtaining closed-form analytic solutions to these equations is impractical (if not unfeasible) in most cases, numerical methods like the \ac{fem} are often employed to compute approximate solutions instead. This method makes use of piecewise polynomial spaces on a partition of the computational domain to approximate the solution of the \ac{pde}. Its many advantages, such as a solid mathematical foundation~\cite{Ern2021}, the capability to handle complex geometries, and the ease of imposing boundary conditions, render \ac{fem} a widely adopted approach for computationally solving \acp{pde}.

Problems posed on complex and/or evolving geometries frequently arise in a wide range of practical applications,
such as those modelled by free-boundary problems.
The use of unstructured mesh generators to generate body-fitted meshes that conform to such geometries is often impractical, challenging, time-consuming, and requiring human intervention in many cases.
Unfitted \acp{fem}, also known as embedded or immersed boundary \acp{fem}, are especially compelling in these scenarios. The differential equation at hand is discretised by embedding the computational domain in an easy-to-generate background mesh (e.g., a simple Cartesian mesh) that does not necessarily conform to its geometrical boundary, thus drastically reducing the geometrical constraints imposed on the meshes to be used for discretisation.
However, unfitted \acp{fem} introduce their own numerical challenges,
such as the small cut cell problem, where arbitrarily
small intersections between the background mesh and the physical domain can lead to numerical instability (i.e., arbitrarily ill-conditioning of the system resulting from discretisation). 
We refer the reader to~\cite{dePrenter2023} for a comprehensive and updated review of unfitted \acp{fem}. 
Two prominent unfitted \acp{fem}, that we leverage in this work, Cut\ac{fem}~\cite{Burman2015} and \ac{agfem}~\cite{Badia2018}, address this issue in different ways: Cut\ac{fem} introduces a ghost penalty term~\cite{Burman2010} to stabilise the system, while \ac{agfem} eliminates ill-posed cells through aggregation techniques and discrete extension operators. Cut\ac{fem} and \ac{agfem} have been successfully applied to a wide range of \acp{pde} on complex geometries, e.g., interface problems~\cite{Neiva2021}, fluid-structure interaction~\cite{Massing2015}, and linear elasticity~\cite{Hansbo2017}.

Recently, \acp{nn} have rapidly gained popularity in the computational science community for approximating solutions to \acp{pde}. The main idea is that one seeks an approximation to the solution of a \ac{pde} from a trial finite-dimensional {\em nonlinear manifold} (e.g., a neural network) as opposed to a finite-dimensional linear space as in \ac{fem}.
\Acp{pinn}~\cite{Raissi2019} are a notable approach for their flexibility and ease of implementation. In this method, the \ac{nn} approximates the solution by minimising the strong \ac{pde} residual at randomly sampled collocation points within the computational domain. An additional penalty term in the loss function evaluated at a set of randomly sampled {\em boundary} collocation points is typically used to weakly impose essential boundary conditions. 
\Acp{pinn} have demonstrated relative success when applied to a variety of \ac{pde}-constrained problems, including fractional \acp{pde}~\cite{Pang2019}, fluid mechanics~\cite{Mao2020}, and inverse problems~\cite{Yang2021}.
Building on a similar concept, \acp{vpinn}~\cite{Kharazmi2021}, which use the weak \ac{pde} residual as the loss function, have been proposed. With $hp$-refinement capability, $hp$-\acp{vpinn} are often more efficient than \acp{pinn} when \ac{pde} solutions feature singularities and steep gradients (see, e.g.,~\cite[Figs.~13 and~14]{Kharazmi2021}). One of the main issues in \acp{vpinn} and \acp{pinn} (apart from the imposition of essential boundary conditions) is the difficulties associated to the accurate integration of \acp{nn} and their derivatives.

To overcome the aforementioned difficulties of \acp{vpinn} and \acp{pinn}, the combination of \acp{nn} with \acp{fe}  has gained significant attention recently. The authors in~\cite{Berrone2022} introduced the \ac{ivpinn} method, which interpolates \acp{nn} onto high-order \ac{fe} spaces to address the inexact integration issues found in \acp{vpinn}. When the same evaluation points for \acp{nn} are used, \acp{ivpinn} showed improved accuracy compared to \acp{vpinn}.
Based on the same interpolation idea, \acp{feinn}~\cite{Badia2024} employ the standard FE offset function to overcome challenges in imposing essential boundary conditions in \ac{pinn}-like methods. Furthermore, \acp{feinn} use a discrete dual norm of the \ac{pde} residual as the loss function, which can be interpreted as a means to precondition the loss to accelerate training convergence and improve accuracy. Similar  techniques have been proposed for \acp{vpinn} in~\cite{Rojas2024}.
\Acp{feinn} can easily adapt to inverse problems by adding a data misfit term into the loss function~\cite{Badia2024}. Another \ac{fe}-based \ac{nn} approach for inverse problems is discussed in~\cite{Xu2025}.
When paired with adaptive mesh refinement, \acp{feinn} can handle \acp{pde} with localised features~\cite{Badia2025adaptive}. With compatible \acp{fe}, they can accurately solve \acp{pde} with weak solutions in $H(\textbf{curl})$ or $H(\textbf{div})$ spaces~\cite{Badia2024compatible}, covering the full range of spaces in the de Rham complex.
Another notable approach combining \acp{nn} and \ac{fe} is the \ac{hidenn}~\cite{Zhang2021}, where \ac{nn} weights and biases are functions of the nodal coordinates in \ac{fe} spaces. This design enables \ac{hidenn} to accomplish $r$-adaptivity during training and $h$-adaptivity by increasing the number of neurons.
The \ac{fenni} framework~\cite{Skardova2024} further improves the efficiency and robustness of \ac{hidenn} training by using a multigrid training strategy. Besides, combining with \ac{pgd}, \ac{fenni}-\ac{pgd}~\cite{DabySeesaram2024} offers a flexible and interpretable solution for real-time surrogate modelling.
In addition to \acp{fem}, \acp{nn} have been combined with other numerical methods, such as finite difference methods~\cite{Xiao2024}.

Several studies in the literature explore the use of \acp{nn} to solve \acp{pde} on complex geometries. PhyGeoNet~\cite{Gao2021} addresses this challenge by mapping the irregular physical domain to a regular reference domain, where classic \acp{cnn} can be trained to approximate the \ac{pde} solution. However, PhyGeoNet lacks the flexibility to handle highly complex domains, as establishing a one-to-one mapping for such geometries is difficult, and direct application of elliptic coordinate transformation is often impossible.
Another popular approach involves training one \ac{nn} in the interior domain and prescribing an offset function that either satisfies or is trained to satisfy the Dirichlet boundary conditions. The final solution is obtained by adding the interior \ac{nn}, multiplied by a distance function that vanishes on the Dirichlet boundary, to the offset function~\cite{Berg2018,Sheng2021}. The shortcoming of this method is that neither the offset function nor the distance function is easily obtainable for complex geometries.
The $\Delta$-\ac{pinn} method~\cite{SahliCostabal2024} enables the use of \acp{pinn} on complex geometries by representing manifolds with eigenfunctions of the Laplace-Beltrami operator. These functions can be approximated numerically for any shape using \acp{fe}. When applied to a steady-state heat transfer problem on a complex sink, the method produces a temperature distribution that is qualitatively closer to the ground truth compared to standard \acp{pinn}. However, the number of eigenfunctions has a significant, non-monotonic impact on the accuracy of the method, and must be carefully tuned for complex geometries.

In this work, we integrate unfitted \acp{fe} into the \ac{feinn} framework to solve \acp{pde} on arbitrarily complex geometries. We coin this novel approach the unfitted \ac{feinn} method.
In a nutshell, \Acp{nn} are interpolated onto high-order unfitted trial \ac{fe} spaces, which are built on a simple background mesh.
The corresponding \emph{linearised} test \ac{fe} space is constructed on a refined mesh. The \acp{nn} are trained to minimise a loss function defined either as the algebraic norm of the residual vector or a discrete dual norm (defined over such test space) of the weak \ac{pde} residual.
The imposition of essential boundary conditions is challenging in \ac{pinn}-like approaches. These methods typically treat the conditions as penalty terms~\cite{Raissi2019}, use Nitsche's method~\cite{Magueresse2024}, or rely on distance functions that vanish at the boundary~\cite{Berg2018,Sheng2021,Sukumar2022}.
Our approach adopts Nitsche's method, as used in unfitted \ac{fem}, and our experiments demonstrate that unfitted \acp{feinn} offer a broad tolerance for the choice of the Nitsche coefficient.
When an Euclidean norm of the discrete residual is used, the method is an extension of the \acp{ivpinn} with Nitsche's method in~\cite{Berrone2023} to unfitted meshes.
However, the use of discrete dual norms proves to be essential for non-smooth solutions. The computation of such norms involves the solution of a Poisson problem discretised using linearised \ac{fe} spaces on unfitted meshes, which is affected by the so-called small cut cell problem~\cite{dePrenter2023}. To solve this issue, both aggregated \ac{fe} spaces~\cite{Badia2018} and ghost penalty stabilisation~\cite{Burman2010} have been considered.

The unfitted \ac{feinn} method transforms the (non)linear system of equations into a residual minimisation problem. The method is insensitive to small cut cell issues when residual vector-based loss functions (that is, based on algebraic norms applied to the residual vector) are employed. Even with loss functions based on the dual norm of the residual functional, only the computation of the discrete dual norm of the residual, i.e., an approximation of the $H^1$ inner product for the linearised \ac{fe} space, requires stabilisation; this is a simpler task compared to stabilising the original unfitted problem with a high-order trial space and a linearised test space.
A comprehensive set of experiments on 2D and 3D complex geometries demonstrates the effectiveness of unfitted \acp{feinn}. We evaluate the method on linear and nonlinear forward \acp{pde}, as well as inverse nonlinear problems, to confirm its reliability for solving \acp{pde} with weak $H^1$ stability.
The geometries considered in the experiments include implicitly defined ones via level-set functions, as well as those defined explicitly by \ac{stl} meshes, which are more representative of real-world applications~\cite{Badia2022b}.
The method is compared against standard \acp{pinn} on both CPU and GPU architectures. Employing the same training points, unfitted \acp{feinn} reach comparable accuracy using fewer training iterations and reduced computational time.
The trained \acp{nn} significantly outperform \ac{fe} interpolations of the analytic solution in terms of $H^1$ error, particularly when the solution is smooth. Furthermore, when a discrete dual norm of the residual is applied, we achieve further gains in accuracy, and the number of iterations required for convergence is substantially reduced.

The rest of the article is organised as follows. \sect{sec:method} introduces the continuous setting, unfitted \acp{fe}, the unfitted \ac{feinn} method, and their associated loss functions. \sect{sec:implementation} details the implementation of linearised test \ac{fe} spaces and provides an overview of the \ac{feinn} implementation. \sect{sec:experiments} presents numerical experiments with different \acp{pde} on 2D and 3D complex geometries. Finally, \sect{sec:conclusions} concludes the article and outlines promising directions for future research.

%% file: method.tex
\section{Methodology} \label{sec:method}
\subsection{Continuous setting} \label{sec:method:continuous}
Let $\Omega \subset \mathbb{R}^d$ be a Lipschitz polyhedral domain with boundary $\partial \Omega$ and outward pointing unit normal $\normal$. We focus on \emph{low-dimensional} problems, i.e., $d\in\{2,3\}$.
The physical domain $\Omega$ may have a complex shape represented implicitly by a level set function $\phi: \mathbb{R}^d \to \mathbb{R}$ such that $ \Omega = \{\pmb{x} \in \mathbb{R}^d: \phi(\pmb{x}) < 0\}$ and $\partial \Omega = \{\pmb{x} \in \mathbb{R}^d: \phi(\pmb{x}) = 0\}$. Alternatively, the boundary $\partial \Omega$ can be explicitly described using \ac{stl} surface meshes~\cite{Badia2022b}.

As a model problem, we consider the Poisson equation with Dirichlet boundary conditions. Its strong form reads: find $u \in U\doteq H^1(\Omega)$ such that
\begin{equation} \label{eq:strong}
  -\pmb{\Delta} u = f \quad \text{in } \Omega, \qquad u = g \quad \text{on } \partial \Omega,
\end{equation}
where $f \in H^{-1}(\Omega)$ is a known source term, and $g \in H^{1/2}(\partial \Omega)$ is a known Dirichlet boundary data. 
Note that Neumann boundary conditions can also be included into the formulation; however, for simplicity and without loss of generality, we focus on the pure Dirichlet case.
The Poisson equation is a representative example of second-order elliptic \ac{pde} with weak $H^1$ stability, which is the class of \acp{pde} considered in this work. Besides, we adopt a nonlinear \ac{pde} from~\cite{Berrone2022} that fits in this class. 
The strong form of the nonlinear \ac{pde} is as follows:
find $u \in U$ such that
\begin{equation} \label{eq:strong_nonlinear}
  -\pmb{\Delta} u + (\pmb{\beta} \cdot \pmb{\nabla}) u + \sigma \mathrm{e}^{-u^2} = f \quad \text{in } \Omega, \qquad u = g \quad \text{on } \partial \Omega,
\end{equation}
where $\pmb{\beta}$ and $\sigma$ are known coefficients.
We emphasise that the methodology in this work can be readily extended to more general \acp{pde}, as far as suitable unfitted \ac{fe} discretisation schemes are available. 

\subsection{Unfitted finite element discretisation} \label{sec:method:fe}
Since $\Omega$ has a complex shape and generating a mesh that fits $\Omega$ is typically challenging, we use unfitted \acp{fe} to discretise the space. We introduce a background box $\Omega^{\rm box}$ that contains $\Omega$ (i.e., $\Omega \subset \Omega^{\rm box}$), and can be easily partitioned by mean of an easy-to-generate mesh, e.g., a simple Cartesian mesh.
Consider a quasi-uniform triangulation $\mathcal{T}_h$ of $\Omega^{\rm box}$ with the characteristic mesh size $h>0$. Following the notation in~\cite{Badia2018}, a cell $K \in \mathcal{T}_h$ can belong to one of three categories: if $K \cap \Omega = \emptyset$, then $K$ is an \emph{external cell} and $K \in \mathcal{T}_h^{\rm ext}$; if $K \cap \Omega = K$, then $K$ is an \emph{interior cell} and $K \in \mathcal{T}_h^{\rm in}$; otherwise, $K$ is a \emph{cut cell} and $K \in \mathcal{T}_h^{\rm cut}$. The set of active cells is defined as $\mathcal{T}_h^{\rm act} \doteq \mathcal{T}_h^{\rm in} \cup \mathcal{T}_h^{\rm cut}$. We then define the standard \ac{fe} space $U_h^{\rm std} \subset U$ of order $k_U$ on $\mathcal{T}_h^{\rm act}$. Note that $U_h^{\rm std}$, a conforming Lagrangian \ac{fe} space on interior and cut cells, is exposed to ill-conditioning due to small cut cells~\cite{Burman2015,Badia2018}.

Unfitted \acp{fe} mainly face three challenges: (i)~numerical integration over the physical domain $\Omega$; (ii)~the imposition of essential boundary conditions; and (iii)~stability issues caused by small cut cells.
We give details on how to address (i) in \sect{sec:implmt:linspace}.
To address (ii), Nitsche's method is widely used to weakly enforce Dirichlet boundary conditions, and we use it in this work as well to build \acp{nn} such that their interpolation onto the trial \ac{fe} space weakly satisfy the essential boundary conditions. 
In the following sections, we first introduce Nitsche's method and then discuss two unfitted \ac{fe} formulations to address (iii).

\subsection{Nitsche's method} \label{sec:method:nitsche}
With unfitted space discretisations, strongly enforcing essential boundary conditions 
in the approximation space (as it is done with body-fitted \acp{fe}) is challenging. Instead, we use Nitsche's method to weakly impose these conditions.
We denote the trial and test \ac{fe} spaces as $U_h$ and $V_h$, respectively. Note that $U_h \subseteq U$ and $V_h \subseteq U$ can either be the standard \ac{fe} spaces introduced in \sect{sec:method:fe} or the aggregated \ac{fe} spaces introduced in \sect{sec:method:unfittedfem}. When necessary, the type of space is indicated with a superscript.
We now define the bilinear form $a: U_h \times V_h \to \mathbb{R}$ and the linear form $\ell: V_h \to \mathbb{R}$ for the Poisson model problem~\eqref{eq:strong} as
\begin{align*}
  a(u_h, v_h) &\doteq \int_{\Omega} \pmb{\nabla} u_h \cdot \pmb{\nabla} v_h \mathrm{d}V + \int_{\partial \Omega} (\frac{\gamma}{h} u_h v_h - u_h \frac{\partial v_h}{\partial \normal} - \frac{\partial u_h}{\partial \normal} v_h) \mathrm{d}S, \\
  l(v_h) &\doteq \int_{\Omega} f v_h \mathrm{d}V + \int_{\partial \Omega} (\frac{\gamma}{h}g v_h - g \frac{\partial v_h}{\partial \normal}) \mathrm{d}S,
\end{align*}
where $\gamma>0$ is the Nitsche coefficient and $\partial/\partial \normal$ denotes derivative in the direction of $\normal$.
These forms are derived by multiplying the Poisson equation~\eqref{eq:strong} by a test function $v_h \in V_h$, integrating by parts while noting that $v_h \neq 0$ on $\partial \Omega$, and adding the Nitsche boundary terms.
Then, the Nitsche \ac{fe} formulation reads: find $u_h \in U_h$ such that 
\begin{equation} \label{eq:nitsche_fe}
  a(u_h, v_h) = \ell(v_h), \quad \forall v_h \in V_h.
\end{equation}
In \ac{fem}, the Nitsche coefficient $\gamma$ in $a$ and $\ell$ is a mesh-dependent parameter that must be sufficiently large to guarantee the coercivity of the bilinear form. Two suggestions to determine $\gamma$ are summarised in~\cite[Sect.~2.1]{Burman2015}.
It is straightforward to verify that Nitsche's method is consistent with the original problem: when $u_h = g$ on $\partial \Omega$, the Nitsche boundary terms in $a$ and $\ell$ vanish, and the formulation reduces to the integration-by-parts form of~\eqref{eq:strong}.
In regards to the convergence of Nitsche's method, we refer to~\cite[Chap.~2]{Thomee2007}.

We denote the weak \ac{fe} residual functional as
\begin{equation} \label{eq:fe_residual}
  \mathcal{R}_h(u_h) \doteq a(u_h, \cdot) - \ell(\cdot) \in V_h'.
\end{equation}

\subsection{CutFEM and AgFEM} \label{sec:method:unfittedfem}
Both Cut\ac{fem} and \ac{agfem} are unfitted methods that use the \ac{fe} formulation~\eqref{eq:nitsche_fe}, the difference lies in the way they stabilise the \ac{fe} scheme. 
Cut\ac{fem}~\cite{Burman2015} utilises the standard \ac{fe} spaces $U_h^{\rm std}$ and $V_h^{\rm std}$ introduced in \sect{sec:method:fe}. To address the small cut cell problem, the Cut\ac{fem} adds ghost penalty terms~\cite{Burman2010} to the left-hand side of~\eqref{eq:nitsche_fe}. 
The penalty is defined on the skeleton $\mathcal{F}_g$, which is the set of interior faces that belong to an element $K$ such that $K \cap \partial \Omega \neq \emptyset$~\cite{Hansbo2017}. Let us denote the normal to a face $F \in \mathcal{F}_g$ as $\normal_F$, the ghost penalty terms are defined as
\begin{equation} \label{eq:ghost_penalty}
  j(u_h, v_h) = \sum_{F \in \mathcal{F}_g} \sum_{i=1}^{k_U} \gamma_{g,i} h^{2i-1} (\jump{\frac{\partial^i u_h}{\partial \normal_F}}, \jump{\frac{\partial^i v_h}{\partial \normal_F}})_F,
\end{equation}
where $u_h \in U_h^{\rm std}$, $v_h \in V_h^{\rm std}$, $k_U$ is the order of $U_h^{\rm std}$, $\gamma_{g,i} > 0$ are the ghost penalty coefficients, $(\cdot, \cdot)_F$ is the $L^2$ inner product over the face $F$, $\partial^i/\partial \normal_F$ is the $i$-th derivative in the direction of $\normal_F$, and $\jump{\cdot}$ denotes the jump in a discontinuous function at $F$:
\begin{equation*}
  \jump{w} \doteq w|_{F_+} - w|_{F_-}.
\end{equation*}
Since high-order directional derivative terms are required in the ghost penalty~\eqref{eq:ghost_penalty}, Cut\ac{fem} becomes computationally expensive when $k_U$ is high. 

Instead of adding additional terms to the formulation, \ac{agfem}~\cite{Badia2018} overcomes this issue at the \ac{fe} space level by cell aggregation and a discrete extension operator. Specifically, ill-posed cut cells are \emph{absorbed} by well-posed interior cells to form aggregates, with each ill-posed cut cell belonging to only one aggregate, and each aggregate containing only one well-posed interior cell, i.e., the so-called root cell of the aggregate. 
The aggregated \ac{fe} spaces $U_h^{\rm ag}$ and $V_h^{\rm ag}$ are then defined on the aggregated mesh $\mathcal{T}_h^{\rm ag}$. In such \ac{fe} spaces, the \acp{dof} in each aggregate are associated only with the root cell, while the other ill-posed \acp{dof} are linear combinations of those in the root cell (the discrete extension from interior to cut cells). We refer to, e.g.,~\cite{Badia2018}, for details of the cell aggregation algorithms, analysis, and implementation details of this method.

\subsection{Neural networks} \label{sec:method:nn}
We focus on fully-connected, feed-forward \acp{nn} composed of affine maps followed by nonlinear activation functions. The network structure is represented by a tuple $(n_0, \ldots n_L)\in \mathbb{N}^{L+1}$, where $L$ is the number of layers and $n_k$ is the number of neurons at layer $1 \leq k \leq L$. We set $n_0 = d$ for the input layer and $n_L = 1$ for the output layer, and $n$ neurons for all hidden layers, i.e., $n_1 = n_2 = ... = n_{L-1} = n$.

At each layer $1 \leq k \leq L$, the affine map is represented by $\pmb{\Theta}_k: \mathbb{R}^{n_{k-1}} \to \mathbb{R}^{n_k}$ and defined as $\pmb{\Theta}_k \pmb{x} = \pmb{W}_k \pmb{x} + \pmb{b}_k$, where $\pmb{W}_k \in \mathbb{R}^{n_k \times n_{k-1}}$ is the weight matrix and $\pmb{b}_k \in \mathbb{R}^{n_k}$ is the bias vector. After the application of the affine map corresponding to each layer except the last one, the activation function $\rho: \mathbb{R} \to \mathbb{R}$ is applied element-wise.
With these definitions, the network is a parametrisable function $\mathcal{N}(\pmb{\theta}): \mathbb{R}^d \to \mathbb{R}$ defined as: 
\begin{equation*}
  \mathcal{N}(\pmb{\theta}) = \pmb{\Theta}_L \circ \rho \circ \pmb{\Theta}_{L-1} \circ \ldots \circ \rho \circ \pmb{\Theta}_1,
\end{equation*}
where $\pmb{\theta}$ represents all the trainable parameters $\pmb{W}_k$ and $\pmb{b}_k$. We apply the same fixed activation function at each layer, even though the activations could be different or trainable.
In this work, $\mathcal{N}$ denotes the \ac{nn} architecture, and $\mathcal{N}(\pmb{\theta})$ represents a particular realisation of this network.

\subsection{Unfitted finite element interpolated neural networks} \label{sec:method:feinn}
The unfitted \ac{feinn} method extends the \ac{feinn} formulation in~\cite{Badia2024} to solve elliptic \acp{pde} on complex domains by using unfitted \ac{fe} formulations~\eqref{eq:nitsche_fe}. We define a \ac{fe} interpolant $\pi_h : \mathcal{C}^0 \rightarrow U_h^{\rm std}$. For the grad-conforming Lagrange (nodal) spaces considered in this work, $\pi_h$ involves only pointwise evaluations at the free nodes. The unfitted \ac{feinn} method aims to find:
\begin{equation} \label{eq:feinn_continuous_loss}
  u_\mathcal{N} \in \argmin{w_\mathcal{N} \in \mathcal{N}}  
  \mathscr{L}(w_\mathcal{N}), \qquad 
  \mathscr{L}(w_\mathcal{N}) \doteq \norm{\mathcal{R}_h({\pi}_h(w_\mathcal{N}))}_{X'},
\end{equation} 
where $u_{\mathcal{N}}$ is the \ac{nn} approximation to the solution of the problem, and $\norm{\cdot}_{X'}$ denotes an appropriate discrete dual norm. 

Let us now discuss the choice of the trial and test \ac{fe} spaces. Since the \ac{fe} problem~\eqref{eq:nitsche_fe} is reformulated as a minimisation problem~\eqref{eq:feinn_continuous_loss}, and the trial \ac{fe} space is ``hidden'' in the \ac{pde} residual~\eqref{eq:fe_residual}, we can choose the trial space to be $U_h^{\rm std}$, i.e., the standard unfitted \ac{fe} space \emph{without any stabilisation}. 
It is not necessary to apply stabilisation because the discrete minimisation problem is determined by the \ac{nn} parameters, not the ill-posed \acp{dof} in the \ac{fe} space.
In fact, all experiments in \sect{sec:experiments} utilise $U_h^{\rm std}$ as the trial space and yield accurate results.
When \acp{nn} are interpolated onto $U_h^{\rm ag}$, they lose control over the cut cells during training because these \acp{dof} are defined as linear combinations of the \acp{dof} in the root cells. 
For the test \ac{fe} space in the residual $\mathcal{R}_h$~\eqref{eq:fe_residual}, the choice depends on the specific loss function employed. Options include $V_h^{\rm std}$, with or without a ghost penalty, or $V_h^{\rm ag}$. Details on the loss functions and their corresponding test space selections are provided in \sect{sec:method:loss}.
Moreover, as shown in~\cite{Badia2024,Badia2025adaptive}, the linearisation of test spaces significantly accelerates convergence and improves the accuracy of \ac{nn} solutions. In unfitted \acp{feinn}, we also employ such Petrov-Galerkin \ac{fe} formulation, adopting \emph{linearised} $V_h^{\rm std}$ and $V_h^{\rm ag}$. These spaces are constructed by refining the background mesh discretising $\Omega^{\rm box}$ and building linear \ac{fe} spaces on the active mesh resulting from intersecting the refined mesh and the domain boundary. Implementation details of these linearised unfitted test spaces are discussed in \sect{sec:implmt:linspace}.

\subsection{Loss functions} \label{sec:method:loss}
The \ac{fe} residual $\mathcal{R}_h$~\eqref{eq:fe_residual} is isomorphic to the residual vector
\begin{equation*}
  [\mathbf{r}_h(w_h)]_i = \left< \mathcal{R}_h(w_h), \varphi^i \right> \doteq \mathcal{R}_h(w_h)(\varphi^i), 
\end{equation*}
where $\{\varphi^i\}_{i=1}^N$ are the shape functions that span the test \ac{fe} space $V_h$ and $N$ is the dimension of $V_h$. We can use the following loss function:
\begin{equation} \label{eq:feinn_discrete_loss}
    \mathscr{L}(u_\mathcal{N}) = \norm{\mathbf{r}_h(\pi_h(u_\mathcal{N}))}_\chi,
\end{equation}
where $\norm{\cdot}_\chi$ is an algebraic norm, e.g., $\ell^1$ or $\ell^2$ norm. When using a body-fitted mesh and the $\ell^2$ norm, the method becomes equivalent to the Nitsche variant of \acp{ivpinn}~\cite{Berrone2023}. However, in contrast to~\cite{Berrone2023}, in this work we focus on unfitted meshes to handle complex geometries. When the loss~\eqref{eq:feinn_discrete_loss} is employed for unfitted \acp{feinn}, the method is referred to as $\ell^1$- or $\ell^2$-\acp{feinn}, depending on the algebraic norm chosen.
Although the $\ell^2$ norm is more common in the literature~\cite{Berrone2022}, the $\ell^1$ norm is also effective and can produce remarkable results, as revealed in~\cite{Badia2024,Badia2025adaptive}. We indeed use both norms in the experiments of this work.
It is important to note that we always choose $V_h^{\rm std}$ as the test space when loss~\eqref{eq:feinn_discrete_loss} is used. In this setup, neither the trial nor test spaces are stabilised, nor is there a need for ghost penalty terms to stabilise the discrete variational \ac{fe} formulation.  Nevertheless, as demonstrated in the moving domain tests in \sect{sec:experiments:2d:moving}, the unfitted \ac{feinn} method remains robust to small cut cells. 
In this setting, the problem is stated as a residual minimisation, where several (local) minima do not disrupt the algorithm. This contrasts sharply with unfitted \ac{fem}, where the solution requires multiplying the inverse of the system matrix and the right-hand side vector, facing ill-conditioning issues, i.e., a lack of uniqueness in the singular limit.

However, as mentioned in~\cite{Badia2024,Badia2025adaptive}, the loss~\eqref{eq:feinn_discrete_loss} becomes ill-posed as $h \rightarrow 0$; it amounts to computing the $L^2$-norm (up to a constant that depends on the quasi-uniformity of the mesh) of the residual functional, which is only in $H^{-1}(\Omega)$. Besides, using these algebraic norms on the residual vector, optimal convergence rates have not been proven~\cite{Berrone2022}. It has motivated the use of dual norms in~\cite{Badia2024}, i.e.,
\begin{equation} \label{eq:feinn_discrete_loss_dual}
  \mathscr{L}(u_\mathcal{N}) = \norm{\mathcal{R}_h(\pi_h(u_\mathcal{N}))}_{U'}.
\end{equation}
In \acp{feinn}, this discrete dual norm can be computed by introducing a discrete Riesz projector $\mathcal{B}_h^{-1}: V_h' \to V_h$ such that
\begin{equation*} 
    \mathcal{B}_h^{-1}\mathcal{R}_h(w_h) \in V_h \ : \left( \mathcal{B}_h^{-1}\mathcal{R}_h(w_h), v_h \right)_U  = \mathcal{R}_h(w_h)(v_h), \quad \forall v_h \in V_h.
\end{equation*}
where $\mathcal{B}_h$ is (an approximation of) the corresponding Gram matrix in $V_h$. Then, we can write the loss \emph{preconditioned} by $\mathcal{B}_h$ as: 
\begin{equation} \label{eq:feinn_precond_loss}
  \mathscr{L}(u_\mathcal{N}) = \norm{\mathcal{B}_h^{-1}\mathcal{R}_h ( {\pi}_h(u_\mathcal{N}))}_X.
\end{equation}
When $X = U$, the loss functions in~\eqref{eq:feinn_discrete_loss_dual} and~\eqref{eq:feinn_precond_loss} are equivalent, and the method is denoted as $H^{-1}$-\acp{feinn}.
Other choices of $X$ are possible, e.g., $X = L^2(\Omega)$. 
To reduce computational costs, $\mathcal{B}_h$ can be substituted with any spectrally equivalent approximation; for example, \ac{gmg} preconditioners are explored in~\cite{Badia2024}.

When utilising the loss function with the discrete dual norm of the residual~\eqref{eq:feinn_precond_loss}, the test \ac{fe} space can be either $V_h^{\rm std}$ or $V_h^{\rm ag}$. However, adding the ghost penalty~\eqref{eq:ghost_penalty} is essential to stabilise the variational problem~\eqref{eq:feinn_precond_loss} in $V_h^{\rm std}$, ensuring its invertibility in inexact arithmetic regardless of the small cuts. In contrast, discretising~\eqref{eq:feinn_precond_loss} in $V_h^{\rm ag}$ is inherently stable~\cite{Badia2018}. Indeed, in~\cite{Verdugo19}, it was shown that standard multigrid solvers for body-fitted \ac{fem} reach expected mesh-independent convergence rates and linear 
complexity when applied to the Gram matrix in $V_h^{\rm ag}$. We compare the performance of these two test spaces in the experiments in \sect{sec:experiments:2d:preconditioners}.

The analysis of \acp{feinn} for the loss function in~\eqref{eq:feinn_discrete_loss_dual} can be found in~\cite[Sect.~3]{Badia2025adaptive} for strong imposition of boundary conditions and coercive systems. The extension to Nitsche's method is straightforward, with minor modifications. These results can be used in the case of unfitted \acp{feinn} under the assumption of a discrete inf-sup condition between the trial and test spaces (see~\cite[Sect.~2.2]{Badia2025adaptive}) combined with the enhanced stability properties of CutFEM and AgFEM discretisation~\cite{Burman2010,Badia2018}.

\subsection{Extension to inverse problems} \label{sec:method:inverse}
Unfitted \acp{feinn} can be extended straightforwardly to solve inverse problems. These include estimating the unknown model parameters of a \ac{pde} from partial and/or noisy observations of the solution and fully reconstructing the solution/state itself.
Let us denote the unknown parameters as $\pmb{\Lambda}$, which may include physical coefficients, boundary values, or source terms. In addition to the state \ac{nn} $u_\mathcal{N}$, we introduce another \acp{nn} $\pmb{\Lambda}_\mathcal{N}$ to approximate the unknown model parameters. 
We interpolate $\pmb{\Lambda}_\mathcal{N}$ onto suitable \ac{fe} spaces using the interpolants $\pmb{\pi}_h$ to ensure efficient computation of the loss function and its gradients. 
The loss function for the inverse problem is defined as: 
\begin{equation} \label{eq:inverse_loss}
  \mathscr{L}(\pmb{\Lambda}_\mathcal{N}, u_\mathcal{N}) \doteq \norm{\mathbf{d} - \mathcal{D}(u_\mathcal{N})}_{\ell^2} + \alpha \norm{\mathcal{R}_h(\pmb{\pi}_h(\pmb{\Lambda}_\mathcal{N}), \pi_h(u_\mathcal{N}))}_{X'}, 
\end{equation}
where $\mathbf{d}\in \mathbb{R}^M$ is the observed data and $M\in\mathbb{N}$ is the number of observations, $\mathcal{D}: U \rightarrow \mathbb{R}^M$ is the observation operator, and $\alpha \in \mathbb{R}^+$ is the penalty coefficient.
The first term in~\eqref{eq:inverse_loss} measures the data misfit, while the second term weakly imposes the \ac{pde} constraint. 
In the inverse problem experiments in \sect{sec:experiments:inverse}, $\mathcal{D}$ is defined by directly evaluating $u_\mathcal{N}$ at the observation points.
We refer to~\cite{Badia2024} for a detailed discussion on \acp{feinn} for inverse problems. The key difference in this work is the use of unfitted \ac{fe} formulations rather than body-fitted \ac{fe} formulations. 

%% file: implmt.tex
\section{Implementation} \label{sec:implementation}
Unfitted \acp{feinn} are a natural extension of the \acp{feinn} framework, thus allowing one to reuse most of the building blocks of the software implementing the framework. The additional components required by unfitted \acp{feinn} are: (i) the design and implementation of unfitted linearised test \ac{fe} spaces and (ii) the integration over the physical domain with irregular shapes when such domain is cut by a background mesh. We first address (i) and (ii) in \sect{sec:implmt:linspace}, followed by an outline of the Julia implementation of unfitted \acp{feinn} in \sect{sec:implmt:feinn}.

\subsection{Unfitted linearised test spaces} \label{sec:implmt:linspace}
The (body-fitted) \ac{feinn} method in~\cite{Badia2024} constructs a trial \ac{fe} space of order $k_U$ on a coarse mesh. It applies $k_U-1$ levels of uniform refinement to this coarse mesh and builds a test \ac{fe} space made of linear 
polynomials on the refined mesh. 
We call such space a \emph{linearised} test \ac{fe} space. It has the same dimension (i.e., number of \acp{dof}) as the trial \ac{fe} space. In practice, such an approach is proven to drastically accelerate training convergence and improve \ac{nn} accuracy~\cite{Badia2024}.

In the unfitted case, we can use (and we indeed use) an analogous approach to construct a linearised test space. In particular, we compute a refined background mesh, and then we intersect this mesh and the physical domain. However, as shown below with an example, this can lead to a refined active triangulation that covers only a strict subregion of the coarse active triangulation. As a consequence, a test space made of linear polynomials on this refined active triangulation may result in a reduced number of \acp{dof} compared to those in the trial space. Fortunately, the residual minimisation approach underlying unfitted \acp{feinn} is flexible enough to naturally accommodate such scenarios.

\begin{figure}[ht]
  \begin{subfigure}{\textwidth}
    \centering
    \begin{tabular}{lll}
      \tikz{\draw[very thick,plum] (0,0) rectangle (1.2em,1.2em);} background box                                             &
      \tikz{\draw[fill=lightblue] (0,0) rectangle (1.2em,1.2em);} coarse active cells                                         &
      \tikz{\fill[black] (0,0) rectangle (0.3em,0.3em);} trial \ac{fe} space \acp{dof}                                          \\
      \tikz{\fill[lightgreen] (0,0) rectangle (1.2em,1.2em);} physical domain                                                 &
      \tikz{\draw[red,dashed,fill=lightred,dash pattern=on 3pt off 1pt]  (0,0) rectangle (1.2em,1.2em);} refined active cells &
      \tikz{\path (0,.5) pic[blue] {cross=2.4pt};} test \ac{fe} space \acp{dof}
    \end{tabular}
  \end{subfigure}

  \centering
  \begin{subfigure}{0.32\textwidth}
    \centering
    \includegraphics[width=0.9\textwidth]{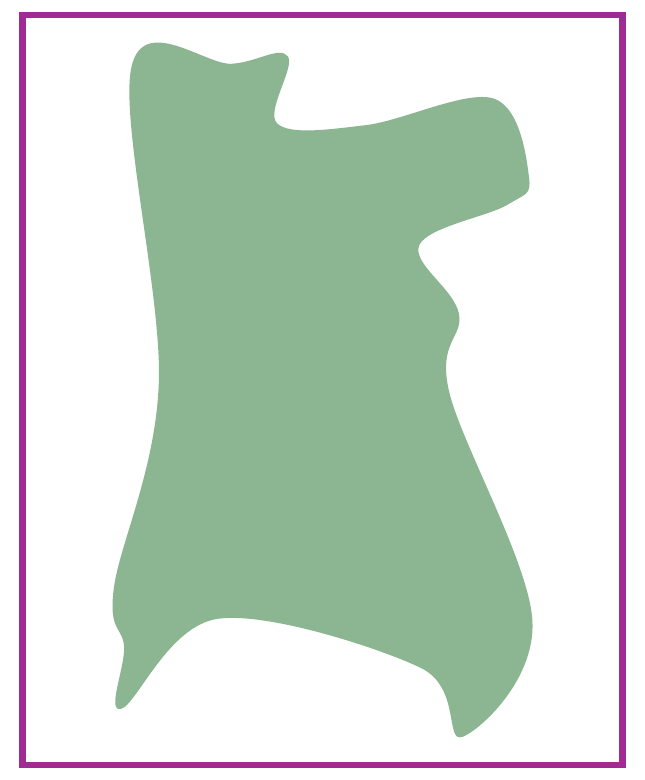}
    \caption{$\Omega^{\rm box}$ and $\Omega$}
    \label{fig:implmt_box_and_phy_domain}
  \end{subfigure}
  \begin{subfigure}{0.32\textwidth}
    \centering
    \includegraphics[width=0.9\textwidth]{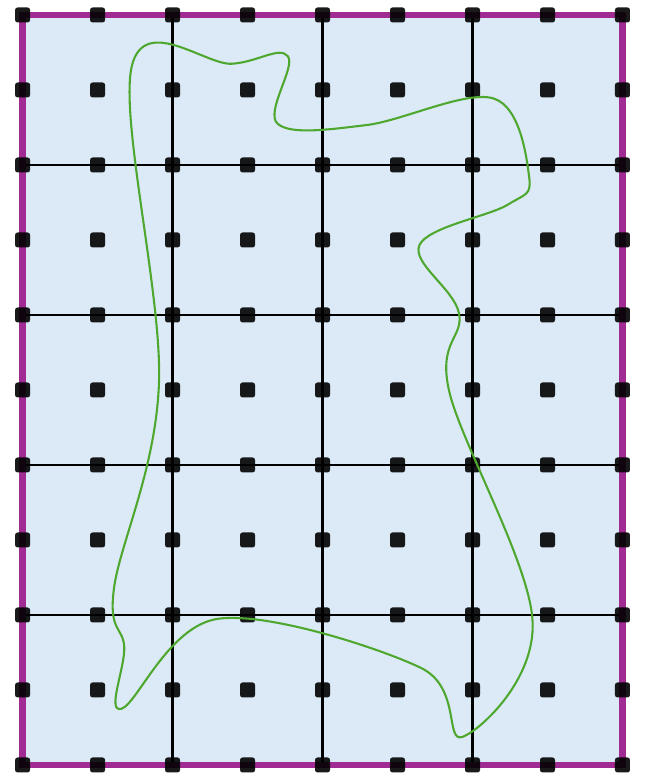}
    \caption{$U_h^{\rm std}$}
    \label{fig:implmt_trial_space_dofs}
  \end{subfigure}
  \begin{subfigure}{0.32\textwidth}
    \centering
    \includegraphics[width=0.9\textwidth]{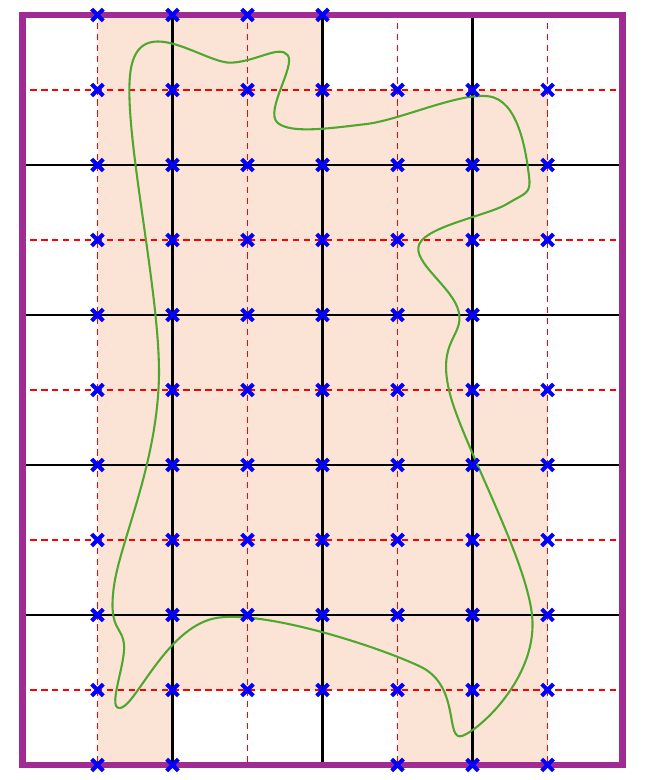}
    \caption{$V_h^{\rm std}$}
    \label{fig:implmt_test_space_dofs}
  \end{subfigure}

  \caption{(a) Background box and physical domain, (b) trial and (c) test \ac{fe} spaces.}
  \label{fig:implmt_domain_and_dofs}
\end{figure}

Consider a physical domain $\Omega$ within a background box $\Omega^{\rm box}$ as displayed in \fig{fig:implmt_box_and_phy_domain}. Let us assume that we discretise $\Omega^{\rm box}$ using a $5\times4$ Cartesian mesh, which we refer to as the coarse mesh. 
As shown in \fig{fig:implmt_trial_space_dofs}, in this particular case, all cells in the coarse mesh are active cells. Assuming $k_U=2$ for the trial space $U_h^{\rm std}$, the global \acp{dof} of $U_h^{\rm std}$ are illustrated by black dots in \fig{fig:implmt_trial_space_dofs}.
Next, we uniformly refine the coarse mesh once, generating a refined mesh with $10\times8$ cells, as shown in \fig{fig:implmt_test_space_dofs}. 
The active cells on this refined mesh, resulting from the  intersection of the refined mesh and the domain,  are highlighted in orange. A linear test \ac{fe} space $V_h^{\rm std}$ is then constructed on this refined active cells, with its \acp{dof} represented by blue crosses in \fig{fig:implmt_test_space_dofs}.

We immediately observe that $V_h^{\rm std}$ has fewer \acp{dof} than $U_h^{\rm std}$ (the same applies to $V_h^{\rm ag}$). A Petrov-Galerkin \ac{fe} discretisation of the problem at hand results in a \emph{rectangular} discrete system. While this mismatch is problematic for \ac{fem}, as it prevents the application of solvers for algebraic systems of equations, we can still apply the unfitted \ac{feinn} method because such method recasts the problem into a residual minimisation problem, which can be solved using a gradient-based optimiser. Although residual minimisation is also available in \ac{fem} (e.g.,~\cite{Guermond2004}), the (body-fitted) \ac{feinn} method can produce more accurate \ac{nn} solutions than \ac{fem} when the solution is smooth, as shown in~\cite{Badia2024,Badia2025adaptive,Badia2024compatible}.
The computation of the discrete dual norm of the residual involves the inversion of a Galerkin discretisation operator 
built using the test space.

\begin{remark}
It is actually possible to build a test \ac{fe} space made of linear polynomials with the same number of \acp{dof} as a higher order trial \ac{fe} space. This requires a suitable \ac{fe} on the coarse mesh active cells. 
In particular, one can use hierarchical Q$_1$-iso-Q$_k$ \acp{fe}~\cite[Sect.~54.4.2]{Ern2021b}. These \acp{fe} 
lead to a well-posed square Petrov-Galerkin discrete system as the basis functions of this \ac{fe} 
have support across the whole coarse cell.
We indeed tested this \ac{fe} and observed results comparable to those achieved by linear bases on the refined active mesh resulting from cutting the refined background mesh with the domain boundary (i.e., the approach that we finally pursue in this paper).
As linear bases on the refined mesh are significantly easier to implement (i.e., one can reuse verbatim existing unfitted \ac{fe} codes), we chose them in this paper.
\end{remark}

In the implementation, we first intersect the refined mesh with $\Omega$ to identify the active refined cells, which allows us to construct $V_h^{\rm std}$. Next, we determine the active cells on the coarse mesh using the following rule: a coarse cell is active if it has at least one active child cell. This way we do not need to compute the 
intersection of the coarse mesh with $\Omega$. The high-order $U_h^{\rm std}$ is then built on these active coarse cells.
The aggregated test \ac{fe} space $V_h^{\rm ag}$, as usual, is constructed from $V_h^{\rm std}$. We classify interior cells as well-posed and cut cells as ill-posed, then form aggregates with one well-posed cell as root.
The space $V_h^{\rm ag}$ is built on these aggregates, with free \acp{dof} determined by the root cells.

Similar to the body-fitted \ac{feinn} method, integration is performed on the \emph{refined} mesh when one uses linearised \ac{fe} spaces. 
For interior cells, the integration process is the same as in fitted \ac{fem}. For a cut cell $K$, we discretise the intersection $K\cap\Omega$ into a local-to-the-cell submesh of triangles (2D) or tetrahedra (3D) and use quadratures rule for simplices with sufficient degree to integrate exactly over these subcells those polynomials defined on the quadrilateral/hexaedral cells. This is the standard approach pursued in unfitted \ac{fem} to integrate polynomials on 
cut cells.

\subsection{Unfitted finite element interpolated neural networks} \label{sec:implmt:feinn}
We reuse the \ac{feinn} implementation as detailed in~\cite[Sect.~5]{Badia2024} for the Euclidean norm of the residual as loss function and inverse problem cases, and in~\cite[Sect.~4.3]{Badia2025adaptive} when using a discrete dual norm of the residual as loss function. For completeness, we highlight the main Julia packages used in the implementation of unfitted \acp{feinn}. 
The \ac{fem} computations are managed by \texttt{Gridap.jl}~\cite{Badia2020, Verdugo2022}, with unfitted extensions handled by \texttt{GridapEmbedded.jl}~\cite{GridapEmbeddedRepo} and \texttt{STLCutters.jl}~\cite{Martorell2024} packages. While \texttt{GridapEmbedded.jl} offers mesh intersection algorithms for implicit geometry representations, such as level set functions, \texttt{STLCutters.jl} handles explicit geometry representations using \ac{stl} meshes. The neural network is implemented using \texttt{Flux.jl}~\cite{Innes2018}. 
We define custom reverse-mode rules in \texttt{ChainRules.jl}~\cite{ChainRules2024} for gradient computation.

%% file: numexp.tex
\section{Numerical experiments} \label{sec:experiments}
We evaluate the accuracy of unfitted \acp{feinn} in solving both linear and nonlinear equations on 2D and 3D complex geometries. The experiments begin with 2D and 3D forward problems, followed by inverse problems in 2D. 
The linear \ac{pde} considered is the Poisson equation with a constant diffusion coefficient, as described in~\eqref{eq:strong}. Its Nitsche \ac{fe} formulation is given in~\eqref{eq:nitsche_fe}.
 Note that the \ac{feinn} framework can readily handle spatially varying coefficients, as demonstrated in~\cite{Badia2024}. In the nonlinear forward and inverse experiments, we focus on the nonlinear equation in strong form~\eqref{eq:strong_nonlinear}.
We first discuss the settings for forward problems.

We employ the same \ac{nn} structure as in~\cite{Badia2024} for all the forward experiments, i.e., $L = 5$ layers, $n=50$ neurons per hidden layer, and $\texttt{tanh}$ activation.
Unless otherwise specified, we use the $\texttt{BFGS}$ optimiser from $\texttt{Optim.jl}$~\cite{Optimjl2018} to train \acp{nn}, following the settings in~\cite{Badia2024,Badia2025adaptive}.
When the  number of \acp{dof} is relatively small (e.g., fewer than 1,000), \acp{nn} can easily overfit, so we reduce the number of training iterations (e.g., to a few hundreds). However, with sufficient number of \acp{dof}, longer training does not negatively affect the training results, so we typically set a large number (e.g., 10,000) of iterations to ensure the \acp{nn} training process has converged. We also note that higher trial space orders require more iterations for the training to converge.
Besides, by default, we use $\ell^2$-\acp{feinn} (i.e., the Euclidean norm of the residual in loss function~\eqref{eq:feinn_discrete_loss}), which has been proven to be effective for smooth solutions~\cite{Badia2024,Badia2025adaptive}. For rough solutions with localised high gradients, we use $H^{-1}$-\acp{feinn}, i.e., loss~\eqref{eq:feinn_precond_loss} with $X=H^1$, which notably improves robustness as shown in \sect{sec:experiments:2d:preconditioners}.

To assess the accuracy of the identified solution $u^{id}$, we compute its $L^2$ and $H^1$ errors as:
\begin{equation*}
  e_{L^2(\Omega)}(u^{id}) = \ltwonorm{u - u^{id}}, \qquad
  e_{H^1(\Omega)}(u^{id}) = \honenorm{u - u^{id}},
\end{equation*}
where $u$ is the true solution, $\ltwonorm{\cdot} = \sqrt{\int_\Omega |\cdot|^2}$, and $\honenorm{\cdot} = \sqrt{\int_\Omega |\cdot|^2 + |\pmb{\nabla}(\cdot)|^2}$.
After training, where the interpolated \ac{nn} (i.e., $\pi_h(u_\mathcal{N})$) is used to compute the loss and gradients, we obtain a \ac{nn} solution (i.e., $u_\mathcal{N}$), $u^{id}$ can be either $\pi_h(u_\mathcal{N})$ or $u_\mathcal{N}$. 
For the $h$- and $p$-refinement error convergence tests, we repeat each experiment 10 times with different initial \ac{nn} parameters to account for the sensitivity of achieved results with respect to \ac{nn} initialisation. Since the variations in the $\pi_h(u_\mathcal{N})$ solutions are negligible, we only show the average errors among these 10 runs in the plots.

To establish a baseline to compare the accuracy of the identified solutions, we interpolate the true solution onto the trial \ac{fe} space and compute the $L^2$ and $H^1$ errors. These interpolations are denoted as $\pi_h(u)$ in the figures and their errors are referred to as $\norm{u - \pi_h(u)}$ in the text.
There are several reasons that prevent us from using \ac{fem} solutions as a baseline for comparison. Firstly, the \ac{fem} solutions are sensitive to the choice of Nitsche coefficient $\gamma$. Secondly, unfitted \acp{fem} often require additional stabilisation techniques, such as ghost penalties in Cut\ac{fem}, which introduce extra parameters that can affect solution accuracy. Lastly, as mentioned in \sect{sec:implmt:linspace}, the Petrov-Galerkin \ac{fe} formulation proposed in this work leads to a rectangular system that cannot be directly solved using standard sparse direct solvers.
In contrast, computing the \ac{fe} interpolation of the true solution is straightforward and unaffected by these hyperparameters. 

For the inverse experiments, we also utilise the \texttt{BFGS} optimiser. As recommended in~\cite{Badia2024}, lighter \ac{nn} structures are adopted, the activation function is set to $\texttt{softplus}$, and the $\ell^1$ norm of the residual vector is employed in the second term of the loss function~\eqref{eq:inverse_loss}. 
We follow the three-step training strategy proposed in~\cite[Sect.~4]{Badia2024} to train the \acp{nn}. In step 1, the state \ac{nn} $u_\mathcal{N}$ is trained using only the data misfit error (first term in~\eqref{eq:inverse_loss}). In step 2, the unknown model parameter \acp{nn} $\pmb{\Lambda}_{\mathcal{N}}$ are trained using the \ac{pde} residual (second term in~\eqref{eq:inverse_loss}) with $u_\mathcal{N}$ fixed. Finally, in step 3, both $u_\mathcal{N}$ and $\pmb{\Lambda}_{\mathcal{N}}$ are jointly trained using the full loss~\eqref{eq:inverse_loss}. 
The first two steps serve as cost-efficient initialisations for the final step, which is the most computationally intensive. Besides, step 3 is divided into several substeps, with the penalty coefficient $\alpha$ gradually increased to ensure the convergence of the training. 
For more details on the training strategy, refer to~\cite{Badia2024} and~\cite[Chap.~17]{Nocedal2006}.

\subsection{2D forward problems} \label{sec:experiments:2d}
\begin{figure}[ht]
  \centering
  \begin{subfigure}{0.33\textwidth}
    \centering
    \includegraphics[width=\textwidth]{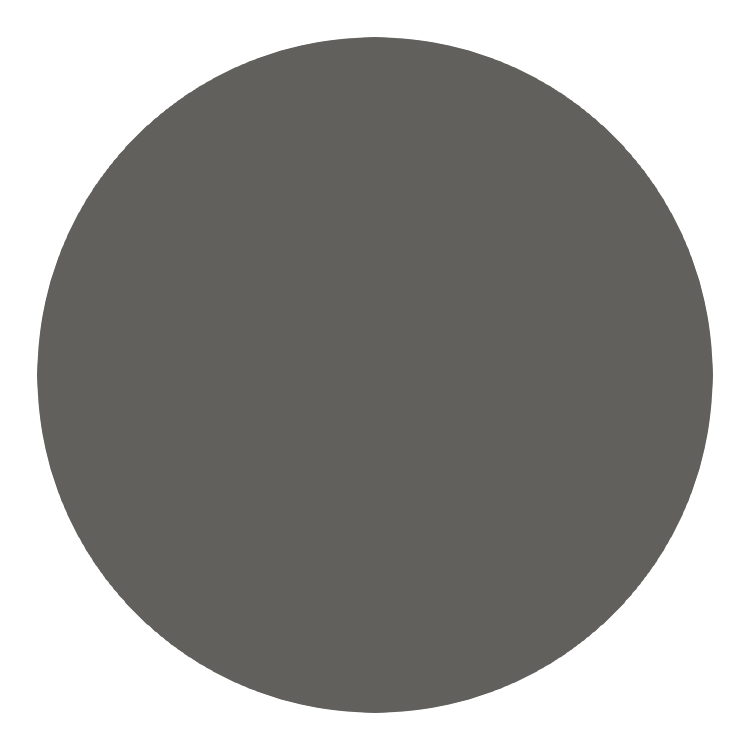}
    \caption{disk}
    \label{fig:geo_disk}
  \end{subfigure}
  \hspace{0.02\textwidth}
  \begin{subfigure}{0.33\textwidth}
    \centering
    \includegraphics[width=\textwidth]{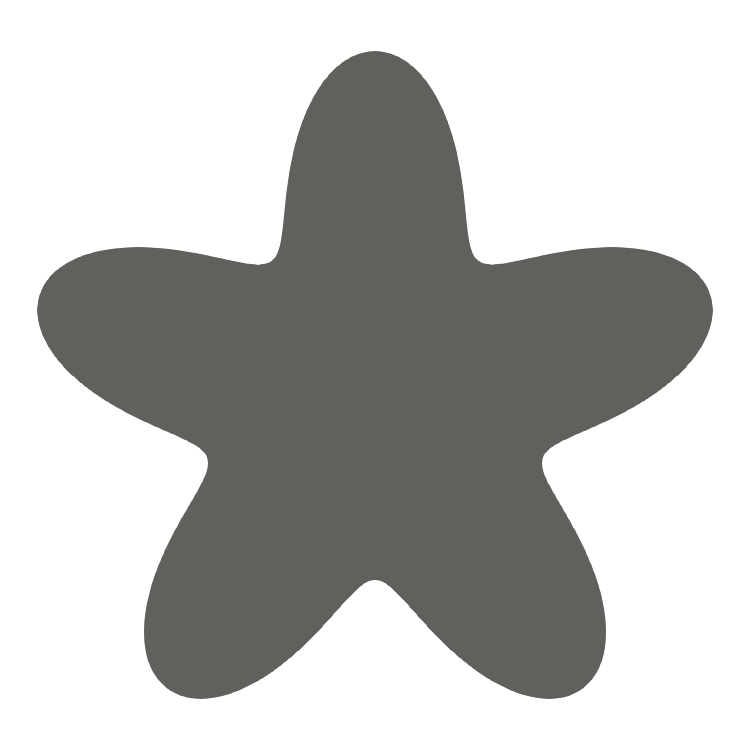}
    \caption{flower}
    \label{fig:geo_flower}
  \end{subfigure}

  \caption{Analytic geometries used in the 2D experiments.}
  \label{fig:geometries_2d}
\end{figure}

As shown in \fig{fig:geometries_2d}, we consider two 2D geometries: a disk and a flower. The disk's level-set function is $\phi(x,y)= \sqrt{(x-x_0)^2 + (y-y_0)^2} - R$, where $(x_0,y_0)$ is the centre and $R$ is the radius. The flower geometry, adapted from~\cite{Neiva2021}, is described by $\phi(\rho, \theta) = \rho - 0.12(\sin(5\theta)+3)$, where $\rho=\sqrt{(x-x_0)^2 + (y-y_0)^2}$ and $\theta$ is the polar angle relative to the centre $(x_0,y_0)$. We use the flower shape in \sect{sec:experiments:2d:nonlinear} and the disk shape in the remaining subsections. 

\subsubsection{The effect of Nitsche coefficient for Poisson equation} \label{sec:experiments:2d:nitsche}
In \ac{fem}, the choice of Nitsche coefficient $\gamma$ significantly affects solution accuracy.
Thus, we first investigate how different $\gamma$ values impact the training of $\ell^2$-\acp{feinn} for solving the Poisson equation.
The physical domain is a disk centred at $(0.5,0.5)$ with a radius of $R=0.4$, and the background mesh is a $50\times50$ uniform grid over $[0,1]^2$. The trial \ac{fe} space has an order of $k_U=3$.
We use different $\gamma$ values in the loss function and train \acp{nn} starting from the same initial parameters for each $\gamma$.
We adopt the following \emph{smooth} solution from~\cite{Badia2024}:
\begin{equation} \label{eq:poisson2d_smooth_solution}
  u(x,y)=\sin(3.2x(x - y))\cos(x + 4.3y) + \sin(4.6(x + 2y))\cos(2.6(y - 2x)).
\end{equation}

\begin{figure}[ht]
  \centering
  \includegraphics[width=\textwidth]{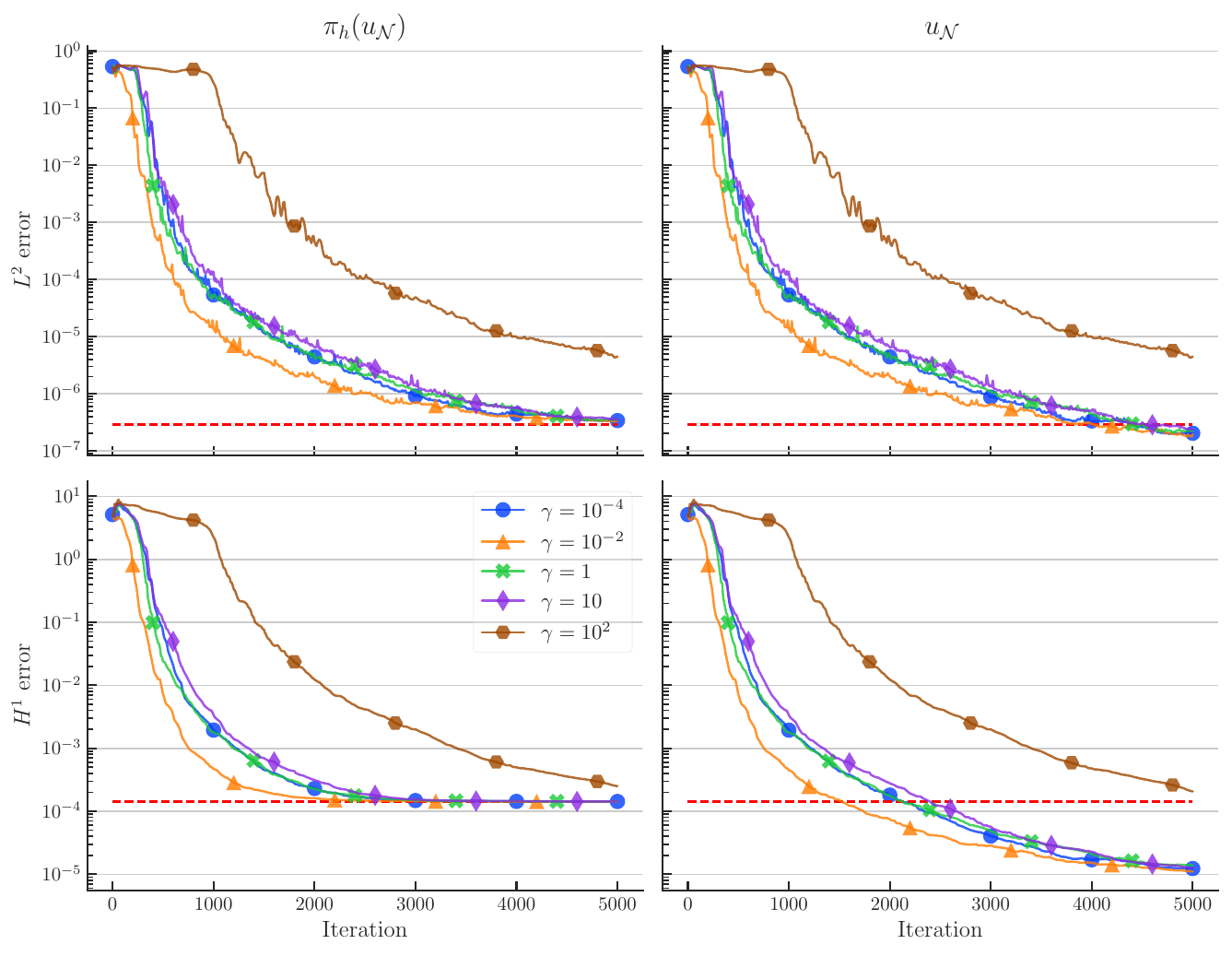}
  \caption{Convergence history of $u^{id}$ in $L^2$ and $H^1$ errors for the 2D Poisson problem, using different Nitsche coefficients during training. The top row shows $L^2$ errors and the bottom row represents $H^1$ errors. The first column displays errors of the interpolated \acp{nn} and the second column corresponds to errors of the \acp{nn} themselves. The red dashed lines represent the \ac{fem} errors.}
  \label{fig:nitsche_error_history}
\end{figure}

\fig{fig:nitsche_error_history} shows the $L^2$ and $H^1$ error histories for the interpolated \acp{nn} and the \acp{nn} themselves during training with different values of $\gamma$ in the loss function.
An immediate observation is that a large $\gamma$ value (e.g., $\gamma=10^2$) leads to a much slower convergence. However, convergence accelerates when $\gamma$ is small, and the corresponding $\pi_h(u_\mathcal{N})$ solutions eventually converge to the $\norm{u - \pi_h(u)}$ lines.
Another important finding is that, consistent with observations in~\cite{Badia2024,Badia2025adaptive}, the \acp{nn} themselves can outperform their \ac{fe} interpolations in terms of $L^2$ and $H^1$ errors. As shown in the second column of \fig{fig:nitsche_error_history}, the curves for $\pi_h(u_\mathcal{N})$ converge to the $\norm{u - \pi_h(u)}$ lines, while the $u_\mathcal{N}$ curves eventually fall below these lines when $\gamma$ is suitable.

Among the tested values, $\gamma=10^{-2}$ achieves the fastest convergence, so for convenience, we simply select $\gamma=10^{-2}$ for the remaining experiments in this work.
However, we emphasise that $\ell^2$-\acp{feinn} are robust to variations in $\gamma$ across a wide range. As shown in \fig{fig:nitsche_error_history}, the differences among convergence curves for $\gamma \leq 10$ values are not significant. For instance, $\gamma=10^{-4}$ and $\gamma=10$ yield similar error histories. Therefore, as long as $\gamma$ is ``small enough'' we expect $\ell^2$-\acp{feinn} to perform well.

\subsubsection{Comparison of \acp{pinn} and unfitted \acp{feinn}} \label{sec:experiments:2d:pinn}
This experiment compares the computational efficiency of unfitted \acp{feinn} with \acp{pinn}. The loss function for \acp{pinn}~\cite{Raissi2019} applied to the Poisson equation~\eqref{eq:strong} reads:
\begin{equation*} \label{eq:pinn_loss}
  \mathscr{L}(u_\mathcal{N}) = \frac{1}{N_{\rm int}} \sum_{i=1}^{N_{\rm int}} \left[\pmb{\Delta}u_\mathcal{N}(x_i, y_i) + f(x_i, y_i)\right]^2 + \tau_b \frac{1}{N_{\rm b}} \sum_{i=1}^{N_{\rm b}} \left[ g(x_i,y_i) - u_\mathcal{N}(x_i,y_i)\right]^2,
\end{equation*} 
where $N_{\rm int}$ and $N_{\rm b}$ are the number of interior and boundary collocation points, respectively, and $\tau_b$ is a penalty coefficient. The first term in the loss imposes the \ac{pde}, while the second term weakly enforces the Dirichlet boundary condition through penalisation.

We adopt the same disk geometry and smooth analytic solution as in \sect{sec:experiments:2d:nitsche}. For the $\ell^2$-\ac{feinn} experiments, the background mesh consists of $25\times25$ quadrilaterals, with $k_U=4$ and $\gamma=10^{-2}$. 
To ensure a fair comparison, we adopt the strategy in~\cite{Berrone2022}, where \acp{pinn} use the same number of collocation points as the number of interpolation nodes employed by \acp{feinn}. In particular, we extract the \ac{feinn} interpolation nodes and use those inside $\Omega$ as interior collocation points for \acp{pinn}. 
Using the level-set function, we then uniformly generate boundary collocation points such that their total number equals to the number the interpolation nodes located outside $\Omega$ in \acp{feinn}. 
In this setup, $N_{\rm int}=5013$ and $N_{\rm b}=868$. This approach ensures both methods are trained on the same interior points and the same amount of total points. However, the computational cost is more favourable to \acp{pinn}, because as opposed to Nitsche's method in \acp{feinn},
they do not require computing spatial derivatives on the boundary (approximately 15\% of the total points).

\begin{figure}[ht]
  \centering
  \includegraphics[width=\textwidth]{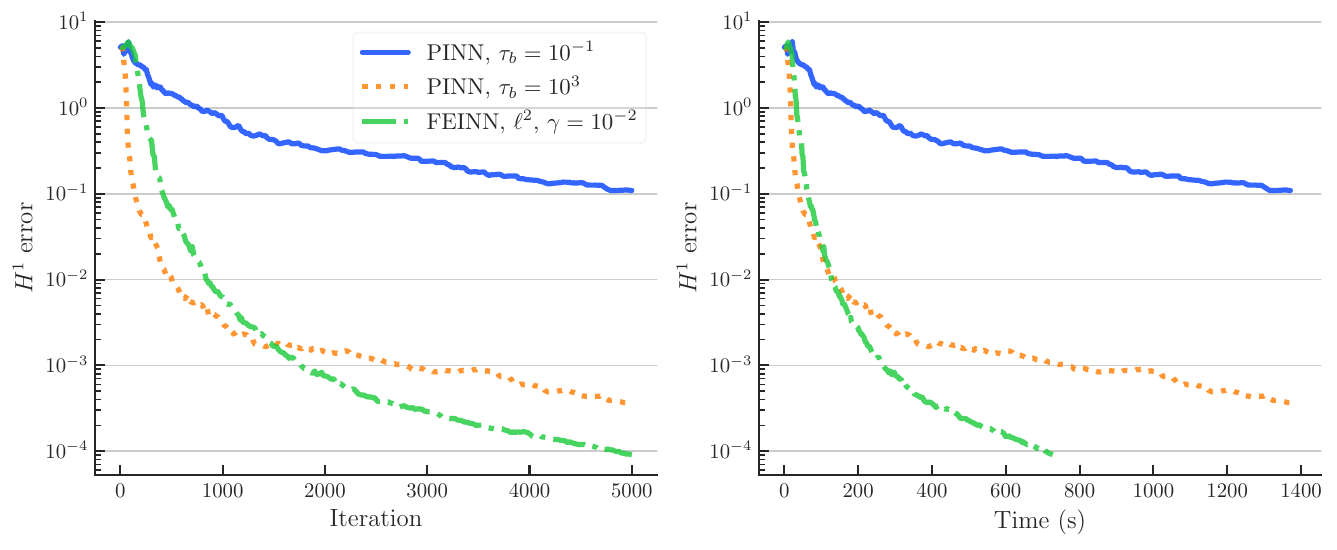}
  \caption{Convergence history of the \ac{nn} $H^1$ error over training iterations and wall-clock time in seconds for the 2D Poisson problem using \acp{pinn} and $\ell^2$-\acp{feinn}. Both methods are implemented in TensorFlow, employ identical \ac{nn} architectures with the same initialisation, and are trained on GPU using TensorFlow's \texttt{BFGS} with an equal number of training points.}
  \label{fig:pinn_cmp}
\end{figure}

As suggested in~\cite{Berrone2022}, convergence in \ac{pinn} training can be accelerated by initially optimising with \texttt{Adam} and then switching to \texttt{BFGS}. 
However, convergence curves in~\cite[Fig. 5]{Berrone2022} show that \texttt{BFGS} yields a steeper descent than \texttt{Adam}. Based on this observation and our preliminary experimental results, we select only \texttt{BFGS} for \acp{pinn} training.
For the penalty coefficient $\tau_b$, we tested several values ranging from $10^{-4}$ to $10^4$ and found that $\tau_b=10^3$ yielded the best results for this problem.
We implement \acp{pinn} in the TensorFlow framework~\cite{tensorflow2015}. 
For a fair comparison, $\ell^2$-\acp{feinn} are also implemented and trained in TensorFlow.\footnote{
We are aware that \acp{pinn} can achieve higher accuracy with alternative optimisers, such as self-scaled BFGS and Broyden methods proposed in~\cite{Urban2025}. We also observe that different implementations of the same optimiser can yield different convergence behaviours. For instance, while \texttt{SciPy}'s \texttt{BFGS} can yield higher accuracy for \acp{pinn} with smaller \ac{nn} architectures, it becomes too inefficient to use for larger \acp{nn} as the ones used in this work. 
A detailed exploration of these optimisers and the subtle differences in their implementations is beyond the scope of this work.}
All experiments are conducted on a single NVIDIA V100 GPU with 32 GB of memory.
The $H^1$ error convergence history of \acp{nn} trained with \acp{pinn} and $\ell^2$-\acp{feinn} is shown in \fig{fig:pinn_cmp}. 

Let us first focus on the left subfigure of \fig{fig:pinn_cmp}, which shows the convergence history of the \ac{nn} $H^1$ error over \texttt{BFGS} iterations.
When the penalty coefficient $\tau_b$ is appropriately chosen (e.g., $\tau_b = 10^3$), \acp{pinn} achieve accuracy comparable to $\ell^2$-\acp{feinn} in terms of $H^1$ error. 
However, the choice of $\tau_b$ significantly affects \acp{pinn}' accuracy: the curve for $\tau_b=10^{-1}$ is placed well above that for $\tau_b=10^3$, indicating that an unsuitable $\tau_b$ can severely hinder \ac{pinn} convergence.
In contrast, as shown in \fig{fig:nitsche_error_history}, \acp{nn} trained by $\ell^2$-\acp{feinn} converge below the $\norm{u - \pi_h(u)}$ reference lines even with suboptimal $\gamma$ values, demonstrating unfitted \acp{feinn}' greater robustness to hyperparameter variation.

The right subfigure of \fig{fig:pinn_cmp} compares the computational time for training \acp{pinn} and $\ell^2$-\acp{feinn}. Notably, the $\ell^2$-\ac{feinn} method completes 5,000 \texttt{BFGS} iterations in around 750 seconds, while \acp{pinn} requires around 1,400 seconds. 
This substantial efficiency gain of $\ell^2$-\acp{feinn} stems not only from avoiding computation of spatial derivatives of \acp{nn} to accelerate loss and gradient evaluations, but also from a better-posed optimisation problem that allows \texttt{BFGS} to converge with fewer gradient calculations.
Specifically, $\ell^2$-\acp{feinn} require 11,308 gradient computations, compared to 13,031 for \acp{pinn} with $\tau_b=10^3$. Moreover, each gradient computation takes 0.0645 seconds for $\ell^2$-\acp{feinn}, versus 0.1064 seconds for \acp{pinn}.
These results align with the findings in~\cite{Berrone2022}, where \acp{ivpinn} were also shown to be faster than \acp{pinn} and \acp{vpinn} in training time.

In summary, for smooth solutions, \acp{pinn} have the potential to achieve similar accuracy as unfitted $\ell^2$-\acp{feinn}; however, their performance is sensitive to the choice of the hyperparameter $\tau_b$. Besides, by avoiding the computation of spatial and nested derivatives of \acp{nn}, $\ell^2$-\acp{feinn} are faster.
It is important to mention that, in this experiment, the level-set function enables easy generation of collocation points for \acp{pinn}. However, for more complex geometries, such as those defined by \ac{stl} meshes, as discussed in \sect{sec:experiments:3d:octopus}, generating collocation points becomes challenging. In contrast, $\ell^2$-\acp{feinn} can handle such complex geometries without additional effort, as demonstrated in the 3D experiments.

\subsubsection{Moving domain test for the Poisson equation} \label{sec:experiments:2d:moving}
We now examine the robustness of unfitted $\ell^2$-\acp{feinn} to changes in the relative position between the physical domain and the background mesh.
We adapt the moving disk experiment from~\cite[Sect.~6.2]{Badia2018}, where a disk with radius $R=0.19$ moves along the anti-diagonal of the background box $[0,1]^2$ by shifting its centre $(x,x)$. As $x$ varies, the relative position of the disk with respect to the background mesh changes, leading to a comprehensive exploration of different cut cell configurations. 
The goal of the experiment is to study how the accuracy of the trained \acp{nn} might be affected by small cut cells.
We use $k_U=3$ and a uniform background mesh with $50\times50$ quadrilaterals. The value of $x$ ranges from 0.2 to 0.8, with a small step size of 0.005 to increase the probability of encountering small cut cells. The true solution is still given by~\eqref{eq:poisson2d_smooth_solution}.

\begin{figure}[htb]
  \centering
  \includegraphics[width=\textwidth]{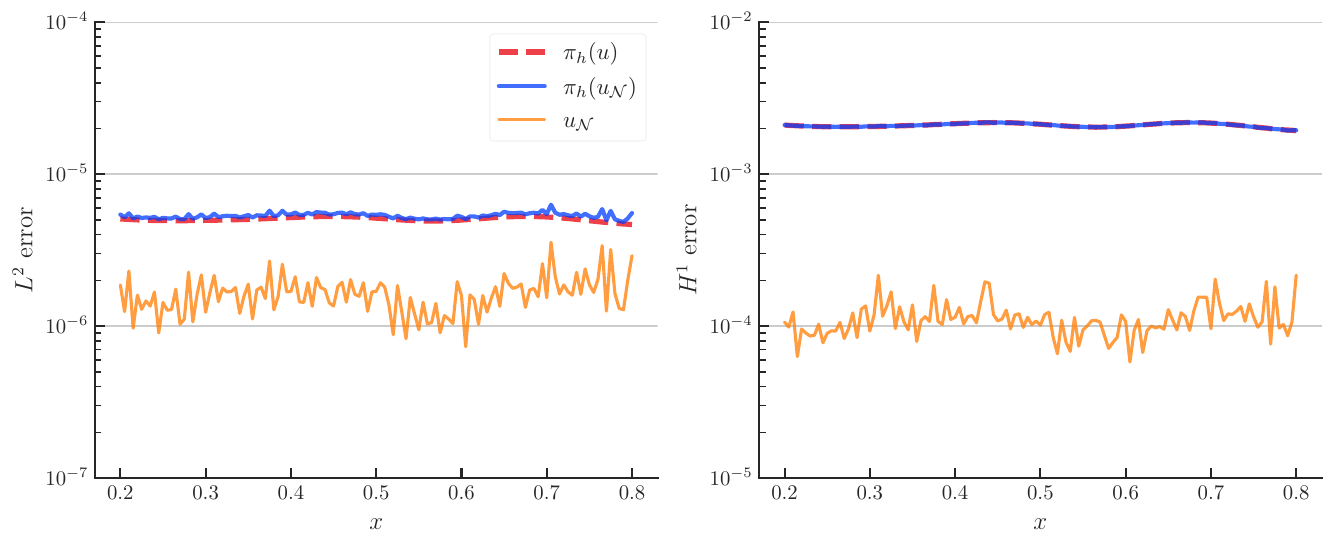}
  \caption{$L^2$ and $H^1$ errors of the trained \acp{nn} for the 2D Poisson problem plotted against the centre position of the moving disk.}
  \label{fig:moving_domain_errors}
\end{figure}

\fig{fig:moving_domain_errors} displays the $L^2$ and $H^1$ errors of the identified solution versus the centre position of the moving disk. The \acp{nn} start with the same initial parameters for each $x$ value.
For the interpolated \acp{nn}, both $L^2$ and $H^1$ errors closely follow the flat $\norm{u - \pi_h(u)}$ lines.
While there are fluctuations in the \ac{nn} errors, they consistently remain below the $\norm{u - \pi_h(u)}$ lines. 
Notably, the $H^1$ errors of the \acp{nn} are steadily more than an order of magnitude lower than the $\pi_h(u_\mathcal{N})$ errors.
These observations highlight that the unfitted $\ell^2$-\ac{feinn} method is highly robust to small cut cells \emph{without using any kind of \ac{fe} stabilisation} such as Ag\ac{fem} or Cut\ac{fem} to compute the residual. In the context of unfitted \acp{feinn}, such stabilisations are required only when using the discrete dual norm in the loss function~\eqref{eq:feinn_precond_loss}, where the discrete dual norms of the residual must be computed.

\subsubsection{Convergence tests for the Poisson equation} \label{sec:experiments:2d:convergence}
In previous experiments, we have established that a Nitsche coefficient of $\gamma=10^{-2}$ yields accurate solutions for $k_U\in \{2,3\}$ and confirmed the robustness of unfitted $\ell^2$-\acp{feinn} against small cut cells. We now confidently proceed to investigate the convergence properties of $\ell^2$-\acp{feinn} with respect to the trial space background mesh size ($h$) and order ($k_U$). We use the disk geometry centred at $(0.5,0.5)$ with radius $R=0.4$ and a uniform Cartesian background mesh on $[0,1]^2$. The true solution is again~\eqref{eq:poisson2d_smooth_solution}.

\begin{figure}[htb]
  \centering
  \includegraphics[width=\textwidth]{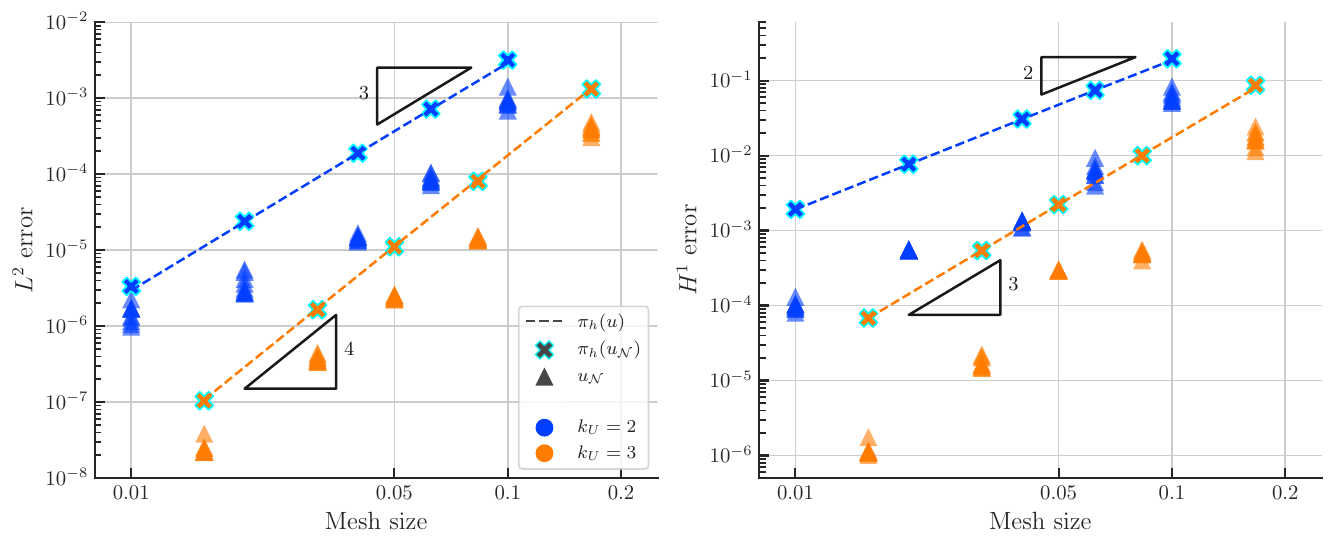}
  \caption{Convergence of the $L^2$ and $H^1$ errors with respect to the trial space mesh size for the 2D Poisson problem. The figure shows results for two orders: $k_U=2$ and $k_U=3$, distinguished by colours.}
  \label{fig:poisson_h_ref_errors}
\end{figure}

We first fix $k_U$ and perform an $h$-refinement experiment for unfitted $\ell^2$-\acp{feinn}.
The convergence results are shown in \fig{fig:poisson_h_ref_errors}. Different colours distinguish between $k_U=2$ and $k_U=3$, while markers are used to differentiate between \acp{nn} and their interpolations.
The dashed lines represent errors of $\pi_h(u)$ versus mesh size.

All the $\norm{u - \pi_h(u)}$ lines in \fig{fig:poisson_h_ref_errors} exhibit the expected slopes, which validates the accuracy of the analytic solution \ac{fe} interpolation and confirms that the error computations use sufficient Gaussian quadrature degree. The $\pi_h(u_\mathcal{N})$ errors always align with the $\norm{u - \pi_h(u)}$ lines for both $k_U=2$ and $k_U=3$, indicating that $\ell^2$-\acp{feinn} converge at the optimal rate.
Both $L^2$ and $H^1$ errors for \acp{nn} are below the corresponding $\norm{u - \pi_h(u)}$ lines, manifesting \acp{nn}' superior accuracy compared to their interpolation counterparts and $\pi_h(u)$. Remarkably, the finest mesh for $k_U=3$ yields $H^1$ errors that are about two orders of magnitude lower than the $\pi_h(u_\mathcal{N})$ errors.
Additionally, the impact of \ac{nn} initialisation is minimal, as the \ac{nn} errors from different initialisations are closely clustered.

\begin{figure}[ht]
  \centering
  \includegraphics[width=\textwidth]{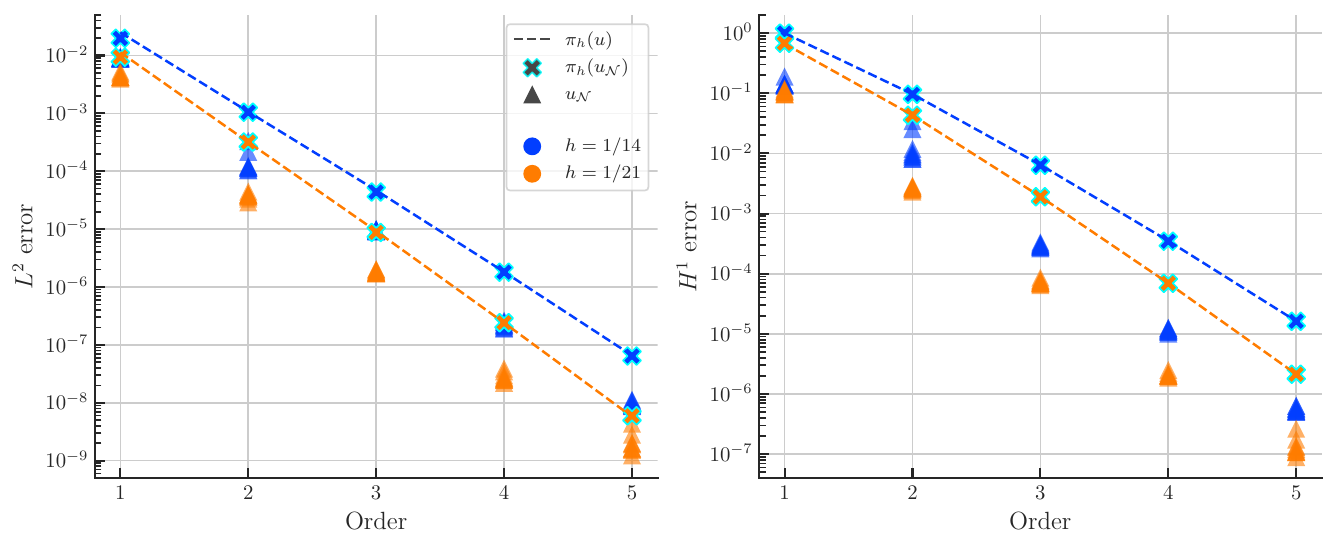}
  \caption{Convergence of the $L^2$ and $H^1$ errors with respect to the trial space order for the 2D Poisson problem. The figure shows results for two mesh sizes: $h=1/14$ and $h=1/21$, distinguished by colours.}
  \label{fig:poisson_p_ref_errors}
\end{figure}

Given the smoothness of the analytic solution~\eqref{eq:poisson2d_smooth_solution}, we also investigate how $k_U$ affects the convergence rate of unfitted $\ell^2$-\acp{feinn} through $p$-refinement tests. We fix the background mesh size and increase $k_U$ from 1 to 5.
The experiment results are displayed in \fig{fig:poisson_p_ref_errors}. We consider two mesh sizes, $h=1/14$ and $h=1/21$, differentiated by colours in the figure.
Similar to the $h$-refinement test, the $\pi_h(u_\mathcal{N})$ errors agree with the true solution interpolation errors, indicating that the method accomplishes the expected convergence rate as $k_U$ increases for both mesh sizes. The \ac{nn} errors are positioned below the $\norm{u - \pi_h(u)}$ lines, demonstrating the improved \ac{nn} accuracy without interpolations. Additionally, \fig{fig:poisson_p_ref_errors} validates the choice of $\gamma=10^{-2}$ for $k_U$ values ranging from 1 to 5, building on the previous subsections where $k_U \in \{2, 3\}$.

In~\cite[Sect.~6.1.1]{Badia2024}, we also conducted $h$- and $p$-refinement tests for the body-fitted \ac{feinn} method on the convection-diffusion-reaction equation, where Dirichlet boundary conditions were strongly imposed via \ac{fe} offset functions. The results showed that \acp{feinn} achieved the expected convergence rates in both tests, and \acp{nn} outperformed \ac{fem} approximations when the analytic solution was smooth.
The current experiments yield similar results for the Poisson equation, further confirming the reliability and accuracy of the \ac{feinn} framework for non-trivial geometries, even when Dirichlet boundary conditions are weakly imposed through Nitsche's method.

\subsubsection{Discrete dual norms for the Poisson equation} \label{sec:experiments:2d:preconditioners}
In this experiment, we explore the impact of discrete dual norms in the training of unfitted \acp{feinn} for the Poisson equation. 
We compare the performance of \acp{pinn} with both unfitted $\ell^2$- and $H^{-1}$-\acp{feinn}. 
We use the same disk geometry as in \sect{sec:experiments:2d:convergence} and set $k_U=3$. The background box $[0,1]^2$ is partitioned into $50\times50$ quadrilaterals.
To increase the difficulty of the problem, we use the following analytic solution, which features a sharp gradient around the centre of the disk:
\begin{equation*}
  u(x,y)=\sin(3.2x(x - y))(5\mathrm{e}^{-100((x - 0.5)^2 + (y - 0.5)^2)} + 1).
\end{equation*}
Since the solution is rougher than~\eqref{eq:poisson2d_smooth_solution}, the use of a discrete dual norm is instrumental in improving the convergence of the training process.

As described in \sect{sec:method:unfittedfem}, we consider two linearised test \ac{fe} spaces: the standard active space $V_h^{\rm{std}}$ consisting of all internal and cut elements, and the aggregated active space $V_h^{\rm{ag}}$ .
The discrete dual norm requires to compute the Gram matrix 
for the linearised test space with respect to the $H^1$-inner product.
However, spectrally equivalent choices can also be considered. For example, in~\cite[Sect.~6.1.3]{Badia2024}, we demonstrated that a single cycle of \ac{gmg} was an effective approximation of the inverse of the Gram matrix.
Since small cut cells in $V_h^{\rm{std}}$ may cause the Gram matrix to become singular in the presence
of rounding errors, additional ghost penalty terms~\eqref{eq:ghost_penalty} are required to ensure the Gram matrix is well-posed.
Note that in the context of unfitted \acp{feinn}, even though the trial space may be of higher order, $V_h^{\rm{std}}$ is always linear. Therefore, a first-order penalty term is sufficient to stabilise the Gram matrix, i.e., $k_U=1$ in~\eqref{eq:ghost_penalty}.
The implementation of unfitted \acp{feinn} does not actually need to deal with high order derivatives in the case of Cut\ac{fem} stabilisation, nor specialised discrete extensions required for high-order aggregated \ac{fe} spaces~\cite{Badia2022c}.
Following the recommendation in~\cite{Burman2015}, we set $\gamma_g=0.1$ in this experiment.
On the other hand, the $V_h^{\rm{ag}}$ space is designed to overcome stability issues by construction, thus no additional penalty terms are needed in the weak formulation.

As shown earlier in \fig{fig:nitsche_error_history}, large Nitsche coefficients can slow down the training convergence, so we include both the suitable $\gamma=10^{-2}$ and a large $\gamma=10^2$ in the experiment to test how effectively $H^{-1}$-\acp{feinn} handles unsuitable Nitsche coefficients.
For clarity in the plots, in this experiment, we only consider the $H^1$ primal norm in the loss~\eqref{eq:feinn_precond_loss}. However, other primal norms can also be considered. For example, in~\cite[Sect.~5.1.2]{Badia2025adaptive}, we showed that four different primal norms were effective for the Poisson equation. 

\begin{figure}[ht]
  \centering
  \includegraphics[width=\textwidth]{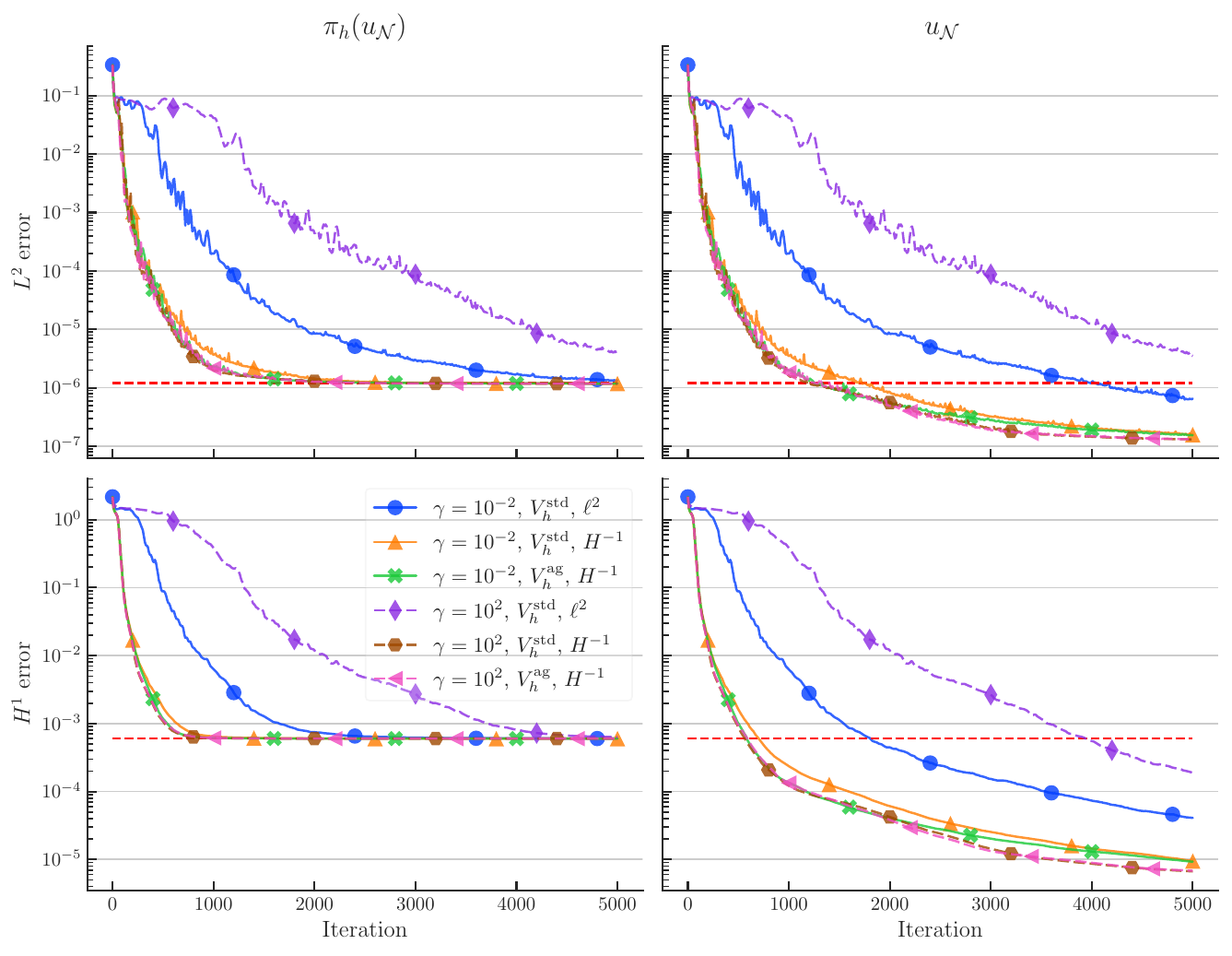}
  \caption{Convergence history of $u^{id}$ in $L^2$ and $H^1$ errors for the 2D Poisson problem, using different combinations of Nitsche coefficient, test space, and norm during training. Refer to the caption of \fig{fig:nitsche_error_history} for details on the information being displayed in this figure.}
  \label{fig:precond_error_history}
\end{figure}

We plot the error histories for $\pi_h(u_\mathcal{N})$ and $u_\mathcal{N}$ for different combinations of test spaces for $\ell^2$-\acp{feinn} and $H^{-1}$-\acp{feinn} in \fig{fig:precond_error_history}.
In the figure, we label $\ell^2$-\acp{feinn} with the postfix ``$V_h^{\rm{std}}$, $\ell^2$'' and $H^{-1}$-\acp{feinn} with the postfixes ``$V_h^{\rm{std}}$, $H^{-1}$'' or ``$V_h^{\rm{ag}}$, $H^{-1}$''. 
We observe that for $\ell^2$-\acp{feinn}, a large Nitsche coefficient also impedes the convergence of both $\pi_h(u_\mathcal{N})$ and $u_\mathcal{N}$ for this problem (with a solution featuring a sharp gradient around the disk centre). However, the use of dual discrete norms in $H^{-1}$-\acp{feinn} overcomes this issue and accelerates the convergence significantly. For instance, with $\gamma=10^2$, the $L^2$ curve of $\pi_h(u_\mathcal{N})$ for $H^{-1}$-\acp{feinn} drops to $10^{-5}$ at around 800 iterations, while $\ell^2$-\acp{feinn} require more than 4,000 iterations to achieve the same error level.
Besides, the use of a discrete dual norm in the loss function is also effective for the suitable Nitsche coefficient ($\gamma=10^{-2}$), as it does not not only accelerates the convergence of \acp{nn} but also enhances the \ac{nn} solution accuracy.
As shown in \fig{fig:precond_error_history}, the corresponding error curves decrease more rapidly than the $\ell^2$-\ac{feinn} ones. In addition, the \ac{nn} error curves for $H^{-1}$-\acp{feinn} in the second column of \fig{fig:precond_error_history} drop below the $\ell^2$-\ac{feinn} ones shortly after training begins.

It is clear from \fig{fig:precond_error_history} that all the dual norm error curves for different $\gamma$ and test space combinations are very close. This implies that dual norm losses are robust to variations in $\gamma$ values and that the choice of test space has only a minor impact on their effectiveness.
Since the Gram matrix requires neither additional stabilisation nor parameter tuning, we recommend using $V_h^{\rm{ag}}$ as the test space in practice.
However, we reiterate that the unfitted \ac{feinn} framework only requires a first-order ghost penalty term to stabilise $V_h^{\rm{std}}$, which makes $V_h^{\rm{std}}$ also a practical choice.

In terms of computational cost, we note that the Gram matrix is symmetric positive definite and is built and factorised using sparse Cholesky factorisation only once before training. Since the problem is linear, this factorisation can be stored and reused throughout training.
At each iteration, the inverse of the Gram matrix is applied to the residual vector by forward and backward substitutions using the sparse Cholesky factor, which makes the training computationally efficient.

Next, we include \acp{pinn} in the comparison to highlight the computational advantages of $H^{-1}$-\acp{feinn}. While their GPU performance was discussed in \sect{sec:experiments:2d:pinn}, we now concentrate on CPU-based training. \Acp{pinn} are still implemented in TensorFlow, whereas both variants of \acp{feinn} are implemented in Julia, as porting our Julia implementation of $H^{-1}$-\acp{feinn} to Python is non-trivial. 
All experiments use identical \ac{nn} architectures with the same initialisation and are run on a single core of an Intel(R) Xeon(R) Platinum 8470Q CPU.
The experimental setup is almost the same as the previous one, except for different initial \ac{nn} parameters.
For $\ell^2$-\acp{feinn}, we use the test space $V_h^{\rm{std}}$ with $\gamma=10^{-2}$. To ensure a fair comparison between $\ell^2$-\acp{feinn} and $H^{-1}$-\acp{feinn}, the later method uses the same test space. Besides, as shown in \fig{fig:precond_error_history}, $H^{-1}$-\acp{feinn} results are very close across different $\gamma$ values. Thus, only the results for $\gamma=10^{-2}$ are presented.

\begin{figure}[ht]
  \centering
  \includegraphics[width=\textwidth]{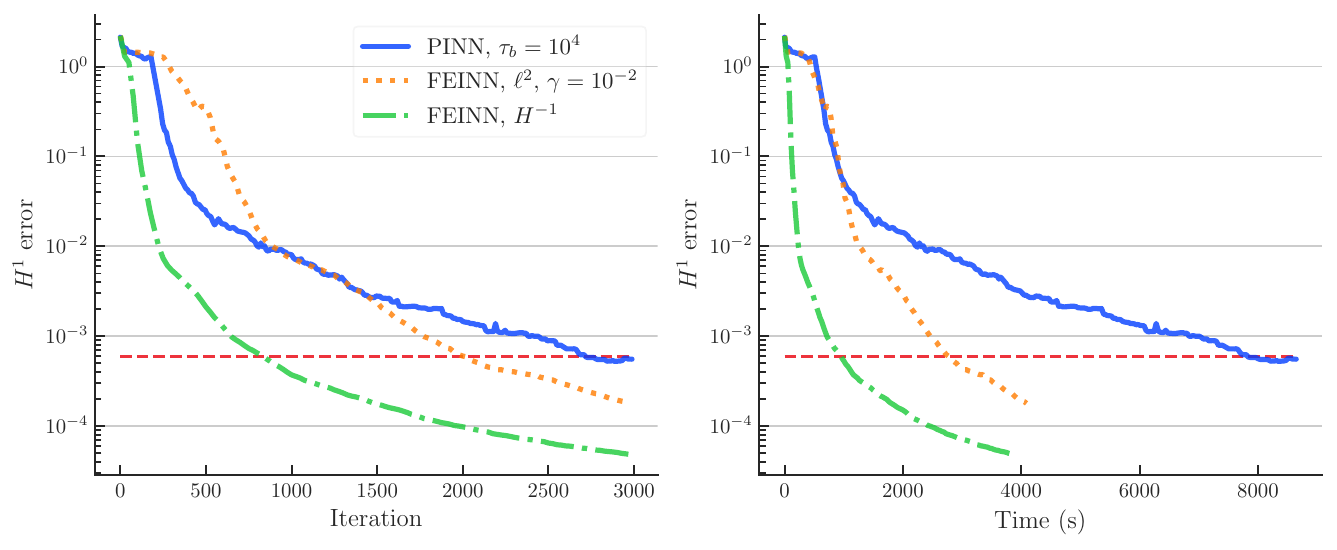}
  \caption{\Ac{nn} error convergence over iteration and time for the 2D Poisson problem, using \acp{pinn}, $\ell^2$-\acp{feinn}, and $H^{-1}$-\acp{feinn}. Both \acp{feinn} variants are implemented in Julia use the standard test space $V_h^{\rm{std}}$. All methods share the same \ac{nn} architecture and initialisation, and are trained on CPU using the \texttt{BFGS} optimiser. The red dashed lines represent the $\norm{u - \pi_h(u)}$ reference lines.}
  \label{fig:comp_time_cmp}
\end{figure}

In \fig{fig:comp_time_cmp}, we present \ac{nn} error history versus training iterations and time in seconds.
Let us first comment on the convergence curves of $\ell^2$-\acp{feinn} and $H^{-1}$-\acp{feinn}.
The use of dual norms indeed improves computational efficiency.
In the $L^2$ error plot on the left of \fig{fig:comp_time_cmp}, even with a suitable Nitsche coefficient ($\gamma=10^{-2}$), $\ell^2$-\acp{feinn} require substantially more training time than $H^{-1}$-\acp{feinn} to reach comparable accuracy. For instance, achieving the reference $\norm{u-\pi_h(u)}$ error, as indicated by the red dashed lines, takes around 2,000 iterations and 3,000 seconds for $\ell^2$-\acp{feinn}, compared to just 800 iterations and 1,000 seconds for $H^{-1}$-\acp{feinn}.
Another notable observation is that the overall computational cost of $\ell^2$- and $H^{-1}$-\acp{feinn} is similar, with both methods finishing 3,000 iterations in around 4,000 seconds. This suggests that incorporating dual norms does not significantly increase the computational cost, while still accelerating convergence. 

In the \ac{pinn} experiments, we again tested different $\tau_b$ values and found that $\tau_b=10^4$ generally yielded the best results.
Among the three methods, \acp{pinn} exhibit the slowest convergence in both time and iterations, requiring more than 2,700 iterations and 7,800 seconds to reach the reference $\norm{u-\pi_h(u)}$ error; $H^{-1}$-\acp{feinn} is noticeably the most robust method. 

\subsubsection{Convergence test for a nonlinear equation} \label{sec:experiments:2d:nonlinear}
As a final forward 2D experiment, we solve the nonlinear equation~\eqref{eq:strong_nonlinear} on the flower geometry (see \fig{fig:geo_flower}), with coefficients $\pmb{\beta} = [2, 3]^{\rm T}$ and $\sigma = 4$. The flower is centred at $(0.5,0.5)$ within the background box $[0,1]^2$.
The true solution is
\begin{equation*}
  u(x,y)=\frac{\cos(\pi(3x + y))}{1 + (x + 2y)^2}.
\end{equation*}
We use the $\ell^1$ norm in the loss~\eqref{eq:feinn_discrete_loss} for this experiment.

\begin{figure}[ht]
  \centering
  \includegraphics[width=\textwidth]{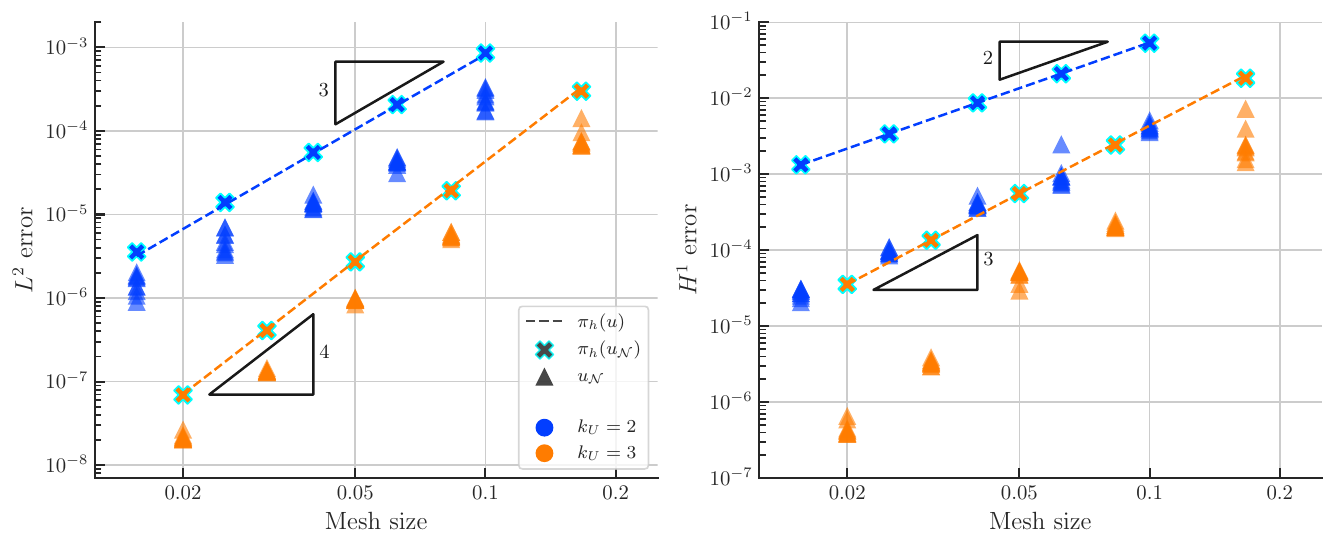}
  \caption{Convergence of the $L^2$ and $H^1$ errors with respect to the trial space mesh size for the 2D nonlinear problem. The figure shows results for two orders: $k_U=2$ and $k_U=3$, distinguished by colours.}
  \label{fig:nonlinear_h_ref_errors}
\end{figure}

Similar to \sect{sec:experiments:2d:convergence}, we conduct $h$-refinement tests for the nonlinear equation on the flower domain. We use $k_U \in \{2,3\}$ and change the background mesh size for each $k_U$. The results are shown in \fig{fig:nonlinear_h_ref_errors}. For this more challenging nonlinear problem, the $\pi_h(u_\mathcal{N})$ errors still exhibit the expected convergence rates, aligning with the $\norm{u - \pi_h(u)}$ lines. The \ac{nn} errors constantly fall below the $\norm{u - \pi_h(u)}$ lines, suggesting a similar accuracy improvement over $\pi_h(u_\mathcal{N})$ as seen in the Poisson equation experiments.

\subsection{3D forward problems} \label{sec:experiments:3d}
\begin{figure}[ht]
  \centering
  \begin{subfigure}{0.33\textwidth}
    \centering
    \includegraphics[width=\textwidth]{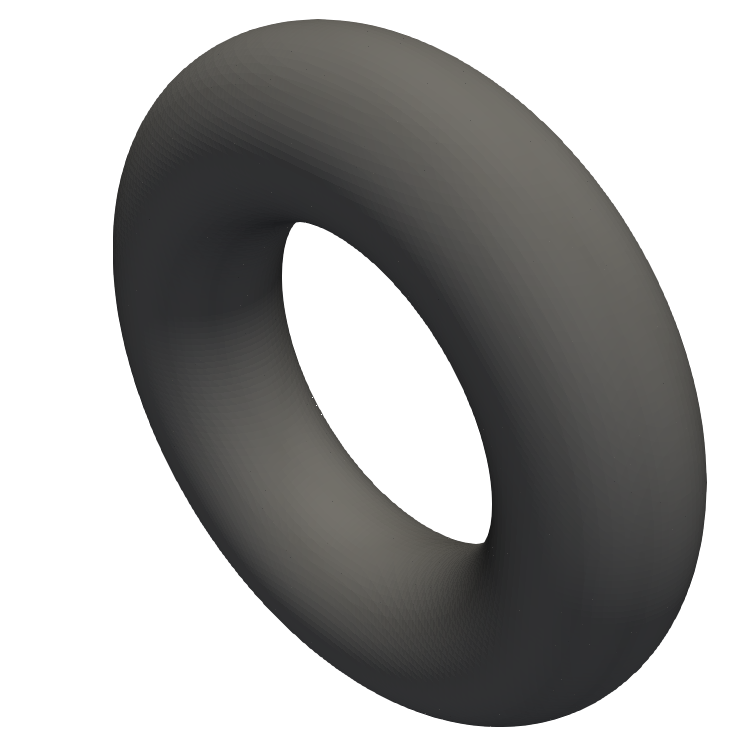}
    \caption{doughnut}
    \label{fig:geo_doughnut}
  \end{subfigure}
  \hspace{0.02\textwidth}
  \begin{subfigure}{0.33\textwidth}
    \centering
    \includegraphics[width=\textwidth]{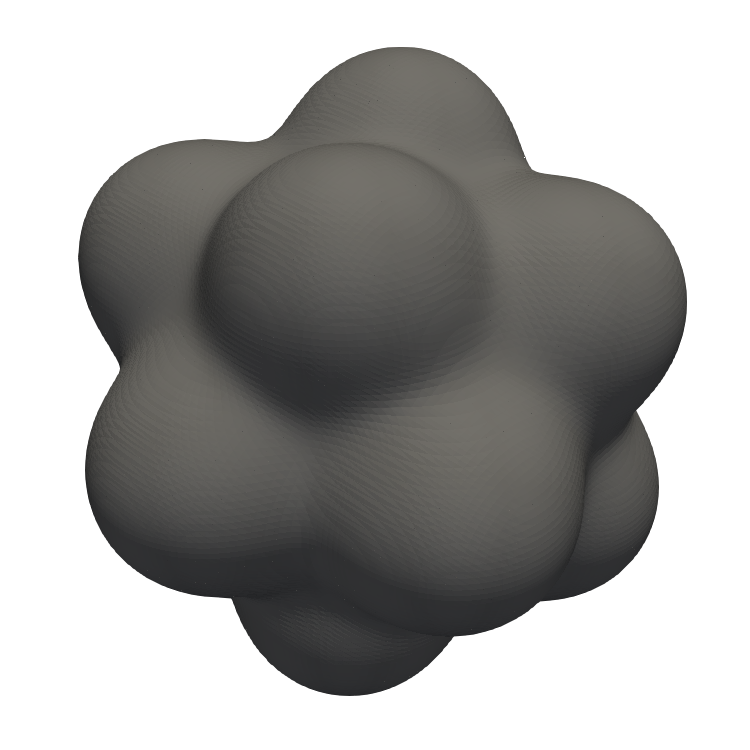}
    \caption{popcorn}
    \label{fig:geo_popcorn}
  \end{subfigure}

  \caption{Analytic geometries used in the 3D experiments.}
  \label{fig:geometries_3d}
\end{figure}

We research two analytic geometries in 3D: a doughnut and a popcorn, as illustrated in \fig{fig:geometries_3d}. Their level-set surface expressions are detailed in~\cite{Burman2015}. For the doughnut, we choose $R=0.3$ and $r=0.1$; for the popcorn, we select $r_0=0.3$, $\sigma=0.09$, and $A=0.9$. Both geometries are centred at $(0.5,0.5,0.5)$.
To further demonstrate the capability of unfitted \acp{feinn} in handling complex 3D geometries, we also include an octopus geometry described explicitly using an \ac{stl} mesh. We apply the $\ell^1$ norm in the loss~\eqref{eq:feinn_discrete_loss} for all 3D experiments.

\subsubsection{Convergence tests on analytic geometries} \label{sec:experiments:3d:convergence}
We first conduct $h$-convergence tests on the doughnut and popcorn geometries. We use trial \ac{fe} space order of $k_U=2$ and alter the mesh resolution for the background box $[0,1]^3$.
The true solution for the Poisson equation is
\begin{equation*}
  u(x,y,z)=\sin(4\pi \sqrt{(x - 1)^2 + (y - 1)^2 + (z - 1)^2}).
\end{equation*}

\begin{figure}[ht]
  \centering
  \includegraphics[width=\textwidth]{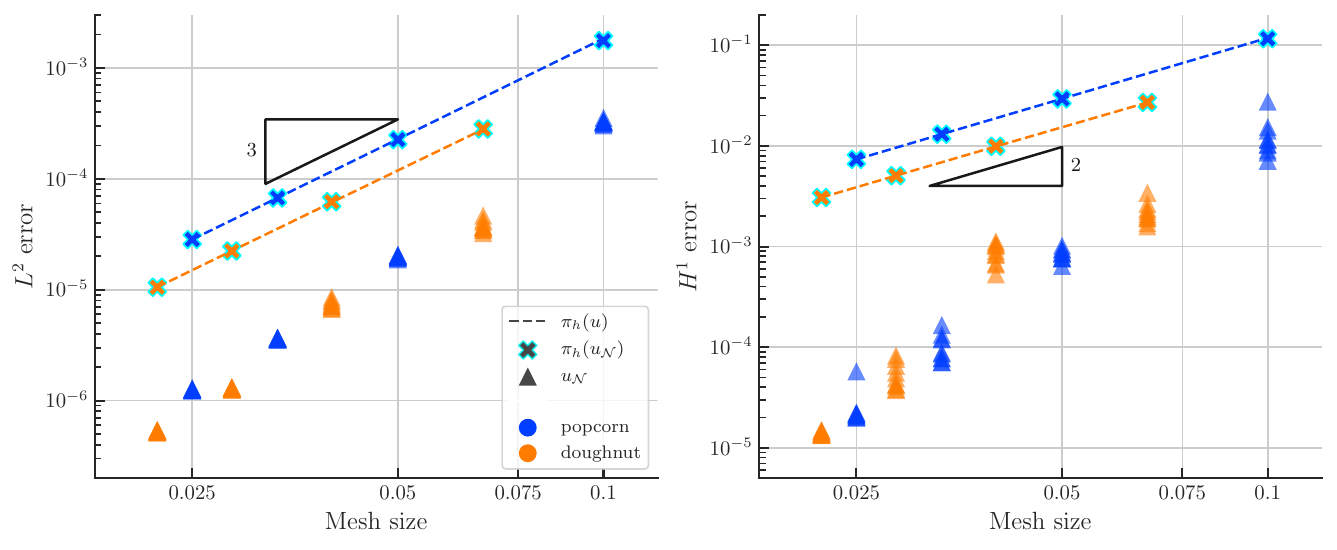}
  \caption{Convergence of the $L^2$ and $H^1$ errors with respect to the trial space ($k_U=2$) mesh size for the 3D Poisson problem. The figure shows results for two geometries: a popcorn and a doughnut, distinguished by colours.}
  \label{fig:poisson3d_h_ref_errors}
\end{figure}

\fig{fig:poisson3d_h_ref_errors} illustrates the decay of $L^2$ and $H^1$ errors for $\pi_h(u_\mathcal{N})$ and $u_\mathcal{N}$ as the mesh size changes, with different colours representing the two geometries. For reference, the $\norm{u - \pi_h(u)}$ dashed lines, showing the expected slopes, are also included.
The $\pi_h(u_\mathcal{N})$ errors clearly achieve the expected convergence rates as they all fall on the $\norm{u - \pi_h(u)}$ lines. The \ac{nn} errors are well below the dashed lines, suggesting superior accuracy of the \acp{nn}, even with complex 3D geometries defined by level-set functions.

\subsubsection{Training on the octopus geometry defined by \ac{stl} meshes} \label{sec:experiments:3d:octopus}
\begin{figure}[ht]
  \centering
  \begin{subfigure}{0.45\textwidth}
    \centering
    \includegraphics[width=\textwidth]{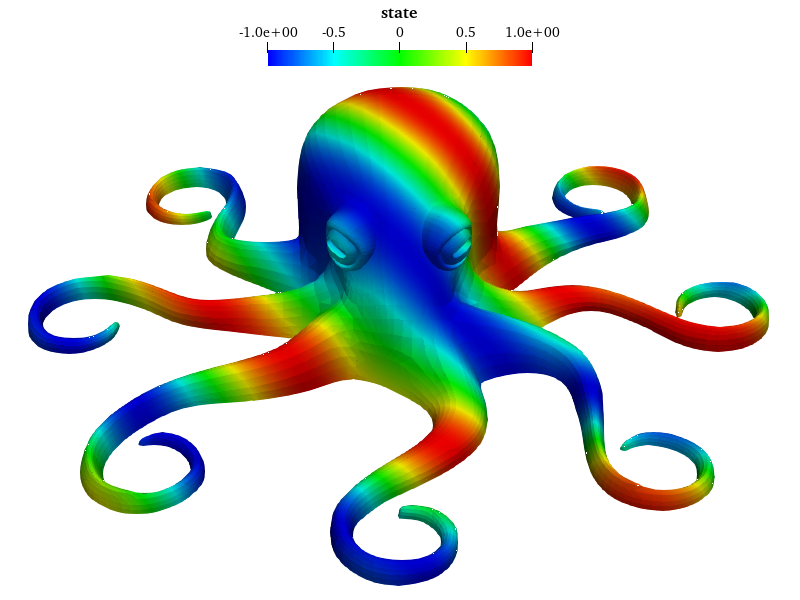}
    \caption{$u$}
    \label{fig:octopus_state}
  \end{subfigure}

  \centering
  \begin{subfigure}{0.45\textwidth}
    \centering
    \includegraphics[width=\textwidth]{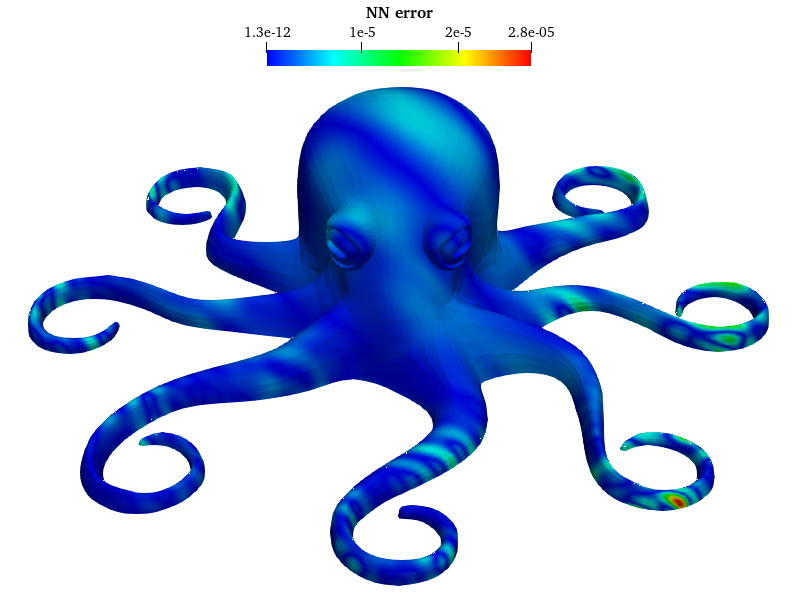}
    \caption{$|u - u_\mathcal{N}|$}
    \label{fig:octopus_nn_error}
  \end{subfigure}
  \hspace{0.01\textwidth}
  \begin{subfigure}{0.45\textwidth}
    \centering
    \includegraphics[width=\textwidth]{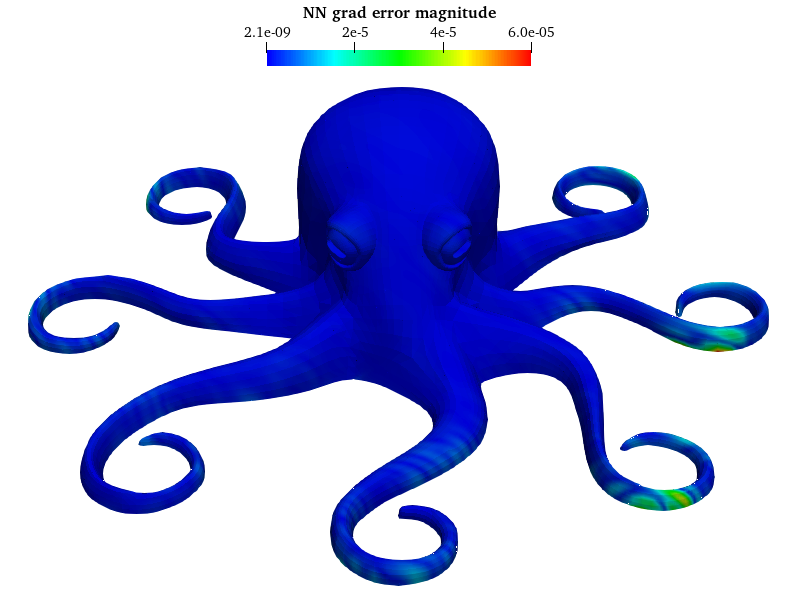}
    \caption{$|\pmb{\nabla}(u - u_\mathcal{N})|$}
    \label{fig:octopus_nn_grad_error}
  \end{subfigure}

  \caption{Reference solution and \ac{nn} errors for the 3D Poisson problem on the octopus geometry.}
  \label{fig:octopus}
\end{figure}

As a concluding experiment for forward problems, we focus on a typical \ac{cae} scenario where the boundary of complex geometries are described by surface \ac{stl} meshes. In practice, the analytic level-set function of a geometry is often unavailable, and the \ac{stl} mesh is a reliable alternative for providing geometric information.
We adopt the octopus geometry, which consists of 10,173 vertices and 20,342 faces, from the Thingi10K dataset~\cite{Zhou2016thingi10k}. This geometry, identified as file ID 69930 in the dataset, is depicted in \fig{fig:octopus}.
We utilise the \texttt{Gridap.jl} ecosystem's Julia package, \texttt{STLCutters.jl}~\cite{Martorell2024}, that provides algorithms to compute the intersection of \ac{stl} meshes and background meshes~\cite{Badia2022b}.
The trial \ac{fe} space of order $k_U=3$ is constructed on a $40\times40\times20$ quadrilateral background mesh over the box $[-80,80]\times[-80,80]\times[-10,40]$.
We solve the Poisson problem, with the true solution illustrated in \fig{fig:octopus_state} and defined by
\begin{equation*}
  u(x,y,z)=\cos(\sqrt[3]{(x - 100)^2 + (y - 100)^2 + (z - 50)^2}).
\end{equation*}

We present the \ac{nn} error in \fig{fig:octopus_nn_error} and its gradient error magnitude in \fig{fig:octopus_nn_grad_error}. The \ac{nn} solution is highly accurate, with point-wise errors below $3\times10^{-5}$ and gradient error magnitudes below $6\times10^{-5}$.
For comparison, \tab{tab:octopus_errors} lists the errors for the \ac{fe} interpolation of the analytic solution $u$, $\pi_h(u_\mathcal{N})$ and $u_\mathcal{N}$ solutions. As expected, the $\pi_h(u_\mathcal{N})$ errors are comparable to the $\pi_h(u)$ errors, while the \ac{nn} errors are notably smaller than those of $\norm{u - \pi_h(u)}$.
\tab{tab:octopus_errors} confirms that unfitted \acp{feinn} are able to accurately approximate the Poisson problem on complex 3D geometries. It also demonstrates that, for smooth solutions, the trained \ac{nn} can achieve higher accuracy than its interpolation, even with highly intricate 3D shapes defined by \ac{stl} meshes. Additionally, it also opens the possibility of using \acp{feinn} for solving problems in real-world \ac{cae} applications.

\begin{table}[ht]
  \aboverulesep=0ex
  \belowrulesep=0ex
  \centering
  \begin{tabular}{c|ccc}
    \toprule 
    & & &\\[-1em]
    & $\pi_h(u)$ & $\pi_h(u_\mathcal{N})$ & $u_\mathcal{N}$\\
    \midrule
    & & &\\[-1em]
    $L^2$ error & $1.78\times10^{-3}$ & $1.51\times10^{-3}$ & $6.33\times10^{-4}$ \\ &&&\\[-1em]
    $H^1$ error & $4.53\times10^{-3}$ & $4.41\times10^{-3}$ & $6.94\times10^{-4}$ \\
    \bottomrule
  \end{tabular}
  \caption{Errors of the \ac{fe} interpolation of the true solution, interpolation of the trained \ac{nn} and \ac{nn} itself for the 3D Poisson problem on the octopus geometry.}
  \label{tab:octopus_errors}
\end{table} 

\subsection{Inverse nonlinear problems} \label{sec:experiments:inverse}
\begin{figure}[ht]
  \centering
  \begin{subfigure}[t]{0.32\textwidth}
    \includegraphics[width=\textwidth]{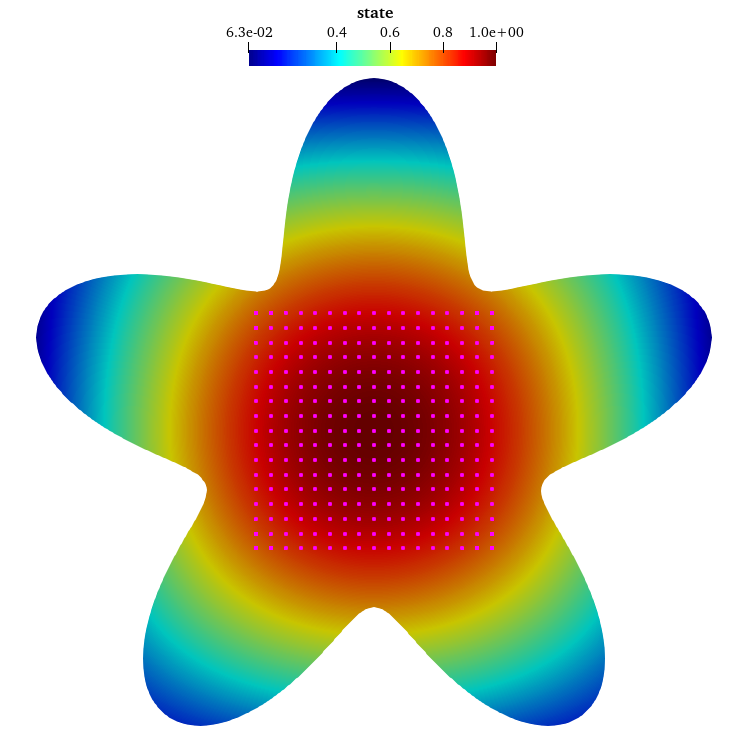}
    \caption{$u$}
    \label{fig:inv_partial_obs_true_state}
  \end{subfigure}
  \begin{subfigure}[t]{0.32\textwidth}
    \includegraphics[width=\textwidth]{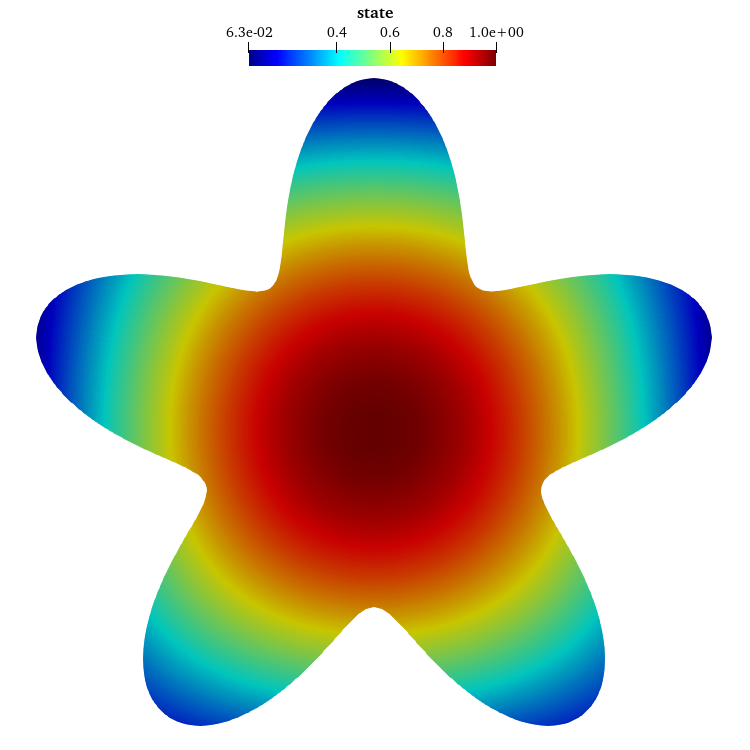}
    \caption{$u_{\mathcal{N}}$}
    \label{fig:inv_partial_obs_nn_state}
  \end{subfigure}
  \begin{subfigure}[t]{0.32\textwidth}
    \includegraphics[width=\textwidth]{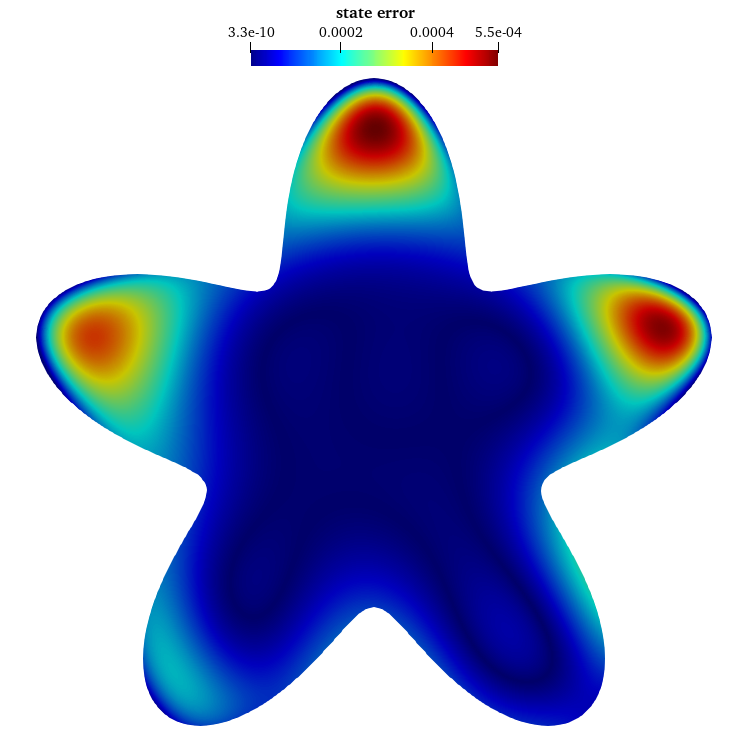}
    \caption{$|u - u_{\mathcal{N}}|$}
    \label{fig:inv_partial_obs_nn_state_error}
  \end{subfigure}

  \centering
  \begin{subfigure}[t]{0.32\textwidth}
      \includegraphics[width=\textwidth]{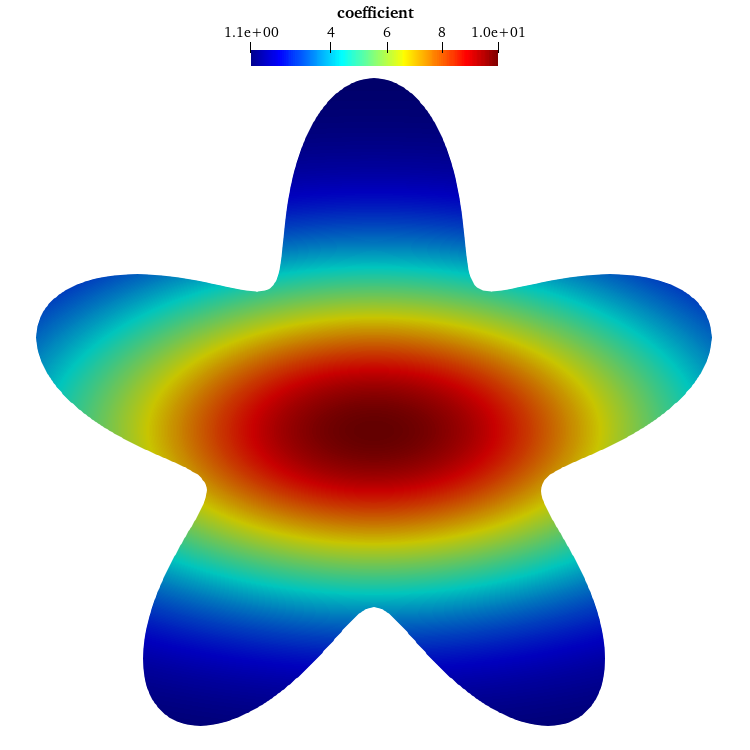}
      \caption{$\sigma$}
      \label{fig:inv_partial_obs_true_coeff}
  \end{subfigure}
  \begin{subfigure}[t]{0.32\textwidth}
    \includegraphics[width=\textwidth]{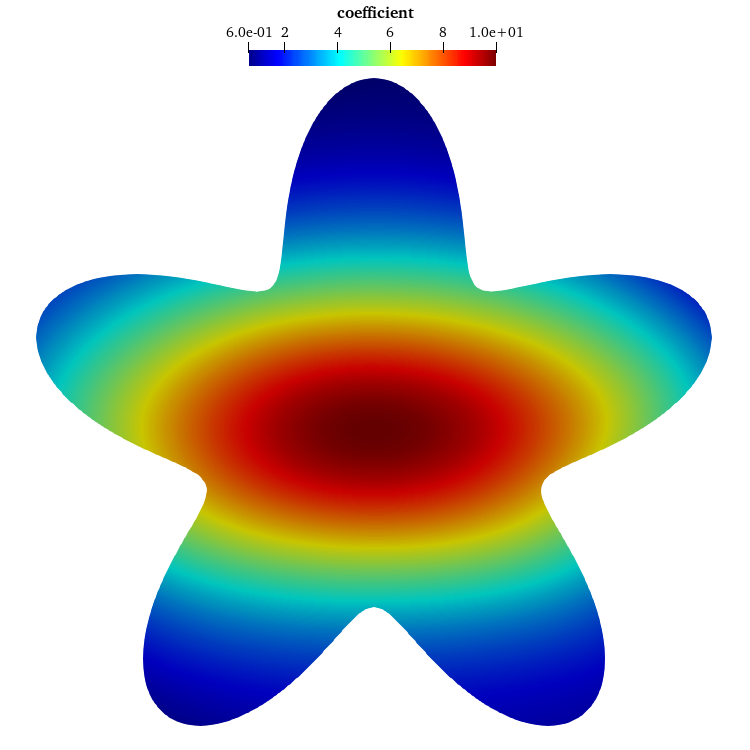}
    \caption{$\sigma_{\mathcal{N}}$}
    \label{fig:inv_partial_obs_nn_coeff}
  \end{subfigure}
  \begin{subfigure}[t]{0.32\textwidth}
    \includegraphics[width=\textwidth]{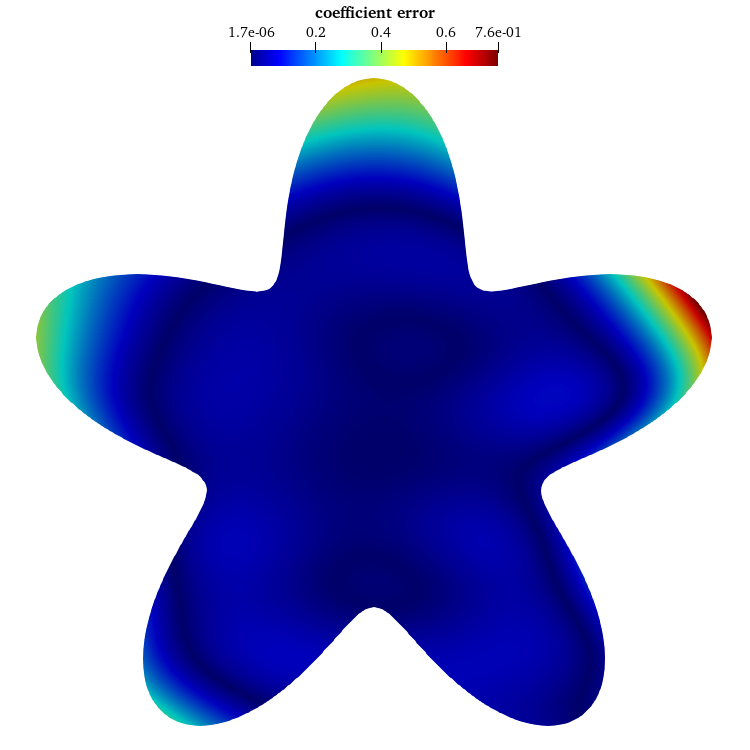}
    \caption{$|\sigma - \sigma_{\mathcal{N}}|$}
    \label{fig:inv_partial_obs_nn_coeff_error}
  \end{subfigure}
    
  \caption{Illustration of analytical (left), \ac{nn} solutions (center) and errors (right) for the inverse nonlinear problem with partial observations. The purple dots in (a) represent the observed data.}
  \label{fig:inv_partial_obs_results}
\end{figure}

We conclude the numerical experiments by solving a nonlinear inverse problem on the flower geometry in \fig{fig:geo_flower}. The objectives are to estimate the coefficient $\sigma$ in the nonlinear equation~\eqref{eq:strong_nonlinear} from partial observations of the state $u$, and to fully reconstruct the partially known $u$. We set $\pmb{\beta}=[2,3]^{\mathrm{T}}$ and define appropriate $f$ and $g$ such that the analytical state and coefficient are given by
\begin{equation*}
  u(x,y)=\sin(\pi x)\sin(\pi y), \quad \sigma(x,y)=1 + 9 \mathrm{e}^{-5((x-0.5)^2 + (2y - 1)^2)}.
\end{equation*}
The first column of \fig{fig:inv_partial_obs_results} displays the true state and coefficient. The observed data are limited to the central region of the domain, indicated by purple dots in \fig{fig:inv_partial_obs_true_state}. Note that we assume prior knowledge of the sign of the unknown coefficient $\sigma$, i.e., $\sigma > 0$.

We adopt the \ac{nn} architectures from~\cite[Sect.~6.2.1]{Badia2024}: both \acp{nn} $u_{\mathcal{N}}$ and $\sigma_{\mathcal{N}}$ have two hidden layers with 20 neurons each; the output layer of $\sigma_{\mathcal{N}}$ is equipped with a rectification function $r(x) = |x| + 0.01$ to ensure positivity. 
All the \ac{fe} spaces, i.e., $u_{\mathcal{N}}$ interpolation or trial space, $\sigma_{\mathcal{N}}$ interpolation space, and test space, are continuous Galerkin spaces of first order on a background mesh consisting of $50\times50$ quadrilaterals. These \ac{fe} spaces are standard unfitted spaces without any stabilisation. Similar to forward problems, essential boundary conditions are weakly imposed via Nitsche's method with $\gamma=10^{-2}$.
The training starts with 400 iterations for $u_{\mathcal{N}}$ in step 1, followed by 100 iterations for $\sigma_{\mathcal{N}}$ in step 2. Step 3 compromises three substeps, with 500 iterations each for the coupled training of $u_{\mathcal{N}}$ and $\sigma_{\mathcal{N}}$. The penalty coefficients $\alpha$ are set to 0.01, 0.03 and 0.09 for these substeps, respectively.

The last two columns of \fig{fig:inv_partial_obs_results} show the \ac{nn} solutions and their point-wise errors. Even though the state observations are limited to partial regions of the domain, the trained $u_{\mathcal{N}}$ in \fig{fig:inv_partial_obs_nn_state} accurately approximates the true state in \fig{fig:inv_partial_obs_true_state} across the entire domain, with point-wise errors below $6\times10^{-4}$ as shown in \fig{fig:inv_partial_obs_nn_state_error}. 
The trained $\sigma_{\mathcal{N}}$ in \fig{fig:inv_partial_obs_nn_coeff} also accurately captures the the true coefficient in \fig{fig:inv_partial_obs_true_coeff}  in most regions of the domain, with point-wise errors in \fig{fig:inv_partial_obs_nn_coeff_error} mostly concentrated around certain parts of the domain boundary. The errors are below $0.2$ in most regions, with the maximum reaching $0.76$ near the boundary.

\begin{figure}[ht]
  \centering
  \includegraphics[width=\textwidth]{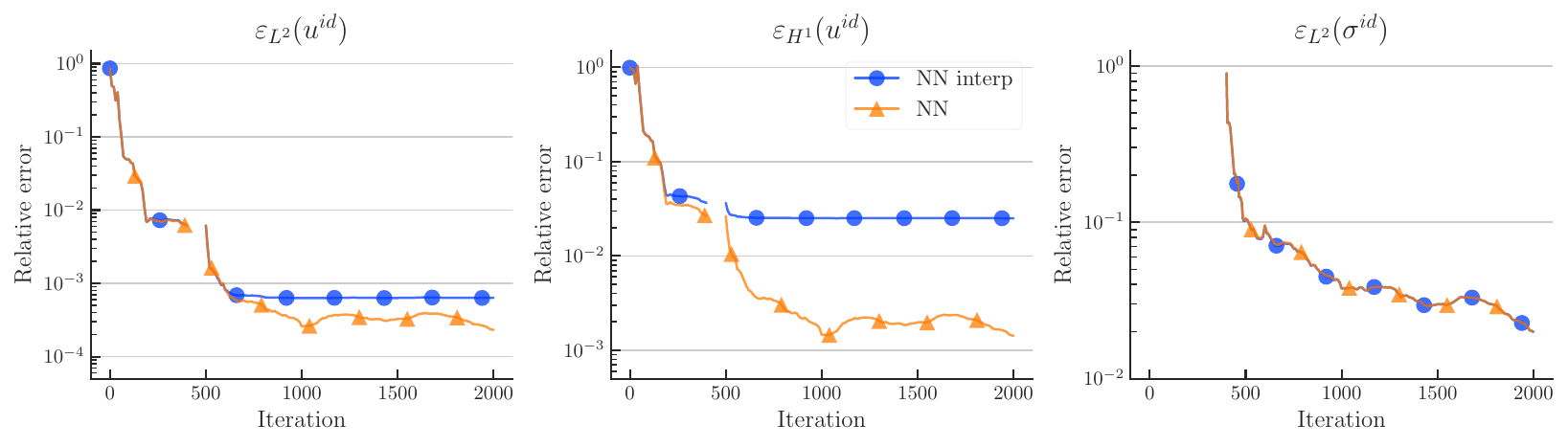}
  \caption{Relative errors during training by the unfitted \ac{feinn} method for the inverse nonlinear problem with partial observations.}
  \label{fig:inv_partial_obs_training_history}
\end{figure}

The relative errors during training for the state and coefficient solutions corresponding to the results in \fig{fig:inv_partial_obs_results} are displayed in \fig{fig:inv_partial_obs_training_history}. The label ``NN interp'' refers to the interpolated \ac{nn} solutions, while ``NN'' refers to the \ac{nn} solutions themselves.
The definition of the relative error is the norm of the difference between the true and identified solutions, divided by the norm of the true solution.
It is evident from \fig{fig:inv_partial_obs_training_history} that the three-step training strategy effectively addresses the inverse nonlinear problem with partial observations. The state solution errors decrease rapidly during step 1 (iterations 1-400), while the coefficient solution errors decrease quickly during step 2 (iterations 401-500). The coupled training in step 3 further reduces the errors in both state and coefficient solutions. 
For the state, the interpolated \ac{nn} solution $\pi_h(u_\mathcal{N})$ stops improving after 700 iterations due to the low-order interpolation of $u_\mathcal{N}$. Despite this, $\pi_h(u_\mathcal{N})$ already provides satisfactory solution, with relative $L^2$ and $H^1$ errors below $0.1\%$ and $3\%$, respectively. Moreover, the state \ac{nn} solution $u_\mathcal{N}$ continues to improve, achieving an order of magnitude gain in accuracy in terms of the relative $H^1$ error.
For the coefficient, both the \ac{nn} solution $\sigma_\mathcal{N}$ and its interpolation $\pi_h(\sigma_\mathcal{N})$ exhibit similar convergence curves, with their relative $L^2$ errors steadily decreasing to $2\%$.

\begin{figure}[ht]
  \centering
  \includegraphics[width=0.65\textwidth]{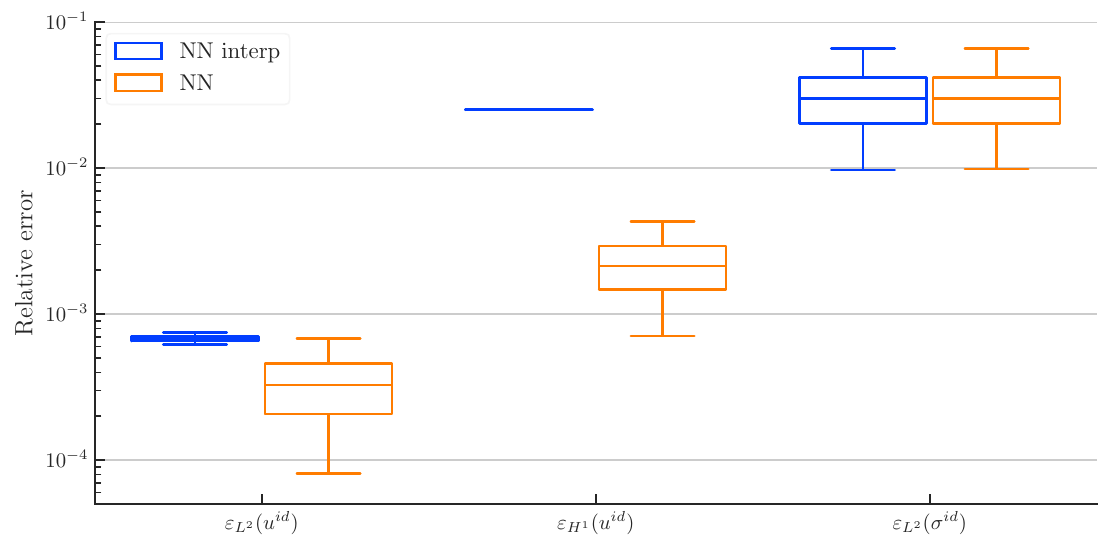}
  \caption{Box plots from 100 experiments for the inverse nonlinear problem with partial observations.}
  \label{fig:inv_partial_obs_boxplots} 
\end{figure}

To further evaluate the robustness of unfitted \acp{feinn} against varying \ac{nn} initialisations, we repeat the experiment 100 times using different initial parameters for both $u_{\mathcal{N}}$ and $\sigma_{\mathcal{N}}$. The box plots of the resulting relative errors are displayed in \fig{fig:inv_partial_obs_boxplots}. 
The interpolated \ac{nn} solutions for the state are consistently accurate across various \ac{nn} initialisations, with almost flat box plots centred at 0.07\% and 3\% for the relative $L^2$ and $H^1$ errors, respectively. 
Although, the \ac{nn} solutions for the state show larger variation compared to the interpolations, their boxes are still relatively compact. Both the relative $L^2$ and $H^1$ error boxes are positioned below those of the interpolated \ac{nn} solutions, indicating improved accuracy. Notably, for the relative $H^1$ error, the third quartile falls below 0.3\%, an order of magnitude lower than the corresponding interpolated \ac{nn} solutions.
The \ac{nn} solutions and their interpolations for the coefficient are nearly identical, with relative $L^2$ error showing first and third quartiles at 2\% and 4\%, respectively.
Overall, \fig{fig:inv_partial_obs_boxplots} confirms the robustness of unfitted \acp{feinn} with respect to different \ac{nn} initialisations in inverse problems with partial observations. This experiment further extends the \ac{feinn} framework from linear inverse problems~\cite[Sect.~6.2]{Badia2024} to nonlinear ones, showcasing the method's high accuracy and reliability on complex geometries.

\section{Conclusions} \label{sec:conclusions}
We combine unfitted \ac{fe} discretisations and the \ac{feinn} framework for the approximation of \acp{pde} on complex geometries using \acp{nn}. In the proposed unfitted \ac{feinn} method, Dirichlet boundary conditions are weakly imposed via Nitsche's method.
The high-order trial \ac{fe} space onto which the \ac{nn} is interpolated is built on the active cells of a coarse background mesh that is unfitted to the boundary of the geometry. With ease of implementation, and existing unfitted \ac{fem} code reuse in mind,  the linearised test space is constructed on the active cells of a uniformly refined version of this coarse background mesh.
Although the trial space has more \acp{dof} than the test space, which leads to a rectangular system, this mismatch 
can be naturally accommodated by the \ac{feinn} method as such framework is based on a residual minimisation principle. Indeed, unfitted \acp{feinn} exhibit outstanding convergence and accuracy in practice.

We conduct extensive numerical experiments to evaluate the unfitted \ac{feinn} method on various 2D and 3D complex geometries for different values of numerical parameters, geometries and \acp{pde}.
As opposed to unfitted \ac{fem}, small cut cells have negligible impact on the accuracy of both trained \acp{nn} and their interpolations, {\em even without using stabilisation techniques} such as Cut\ac{fem} or Ag\ac{fem} to compute the residual (such stabilisations are required, though, for the computation of discrete dual norms). 
The $h$- and $p$-convergence tests for the Poisson equation manifest that the trained \acp{nn} outperform the \ac{fe} interpolation of the analytic solution on the same trial space, and their interpolations accomplish the expected convergence rates. Besides, unfitted \acp{feinn} demonstrate robustness to the choice of Nitsche coefficient, which is a another advantage over unfitted \acp{fem}.

In experiments comparing unfitted \acp{feinn} and \acp{pinn}, \acp{feinn} achieve comparable or better accuracy, demonstrate higher robustness to hyperparameter choices, and require significantly less training time on both GPUs and CPUs. 
We observe that the use of a discrete dual norm of the residual in the loss function, combined with a robust unfitted discretisation, significantly accelerates \ac{nn} training convergence compared to the Euclidean norm variant and \acp{pinn}. This effect is particularly notable for Poisson problems featuring a sharp gradient in a localised region of the domain, substantially reducing computational cost.
The nonlinear equation experiments verify that unfitted \acp{feinn} are capable of dealing with more challenging \acp{pde}.
The experiments on 3D geometries show that unfitted \acp{feinn} accurately solve \acp{pde} on 3D complex geometries defined implicitly via level-set functions or explicitly through \ac{stl} surface meshes. This makes the method well-suited for real-world \ac{cae} scenarios. 

The inverse problem experiments exhibit that unfitted \acp{feinn} can be easily and effectively extended to solve nonlinear inverse problems, delivering highly accurate and robust solutions on complex geometries. This overcomes key limitations of traditional \acp{fem} for inverse settings, which require the use of adjoint solvers, nested loops, and regularisation.

Future research could explore several promising directions. 
First, as discussed in~\cite{Badia2025adaptive}, to exploit the nonlinear approximability properties of \acp{nn}, it is essential to adaptively refine the \ac{fe} spaces. Enhancing unfitted \acp{feinn} with adaptive training techniques~\cite{Badia2025adaptive} will enable us to tackle problems with multiscale features or singularities on complex geometries.
Another encouraging research avenue is applying unfitted \acp{feinn} to evolutionary \acp{pde}, especially time-dependent problems on moving domains~\cite{Badia2023} and fluid-structure interaction problems with evolving geometries~\cite{Massing2015}. In time-dependent problems, the \ac{nn} trained at the previous time step serves as a good initial guess for the current step, which can clearly improve the training efficiency of unfitted \acp{feinn}.